\documentclass[12pt]{article}
\usepackage{latexsym, amssymb, amsmath, amscd, amsfonts, epsfig, graphicx, colordvi,verbatim,ifpdf,tikz,extarrows}
\usepackage{amsfonts, amsmath, amssymb}
\usepackage{amssymb,amsfonts,amsmath,latexsym,epsfig,cite, psfrag,eepic,color}
\usepackage{amscd,graphics}
\usepackage{latexsym, amssymb,  amsmath,amscd, amsfonts, epsfig, graphicx, colordvi,amsthm}

\usepackage[colorlinks,linkcolor=blue,anchorcolor=blue,citecolor=blue]{hyperref}

\usepackage{graphicx}
\usepackage{color}
\usepackage{ifpdf}
\usepackage{fancybox}

\usepackage{float}

\newtheorem{thm}{Theorem}[section]

\newtheorem{conj}[thm]{Conjecture}

\newtheorem{lem}[thm]{Lemma}
\newtheorem{core}[thm]{Corollary}

\def\pf{\noindent{\it Proof.} }
\def\qed{\nopagebreak\hfill{\rule{4pt}{7pt}}
\medbreak}

\hypersetup{colorlinks = true}

\newcommand{\ospt}{\mathop{\mathrm{ospt}}\nolimits}

\numberwithin{equation}{section}

\def\qed{\nopagebreak\hfill{\rule{4pt}{7pt}}
\medbreak}

\newlength{\boxedparwidth}
  {\begin{center} \begin{tabular}{|@{\hspace{.315in}}c@{\hspace{.15in}}|}
                  \hline \\ \begin{minipage}[t]{\boxedparwidth}
                  \setlength{\parindent}{.25in}}%
  {\end{minipage} \\ \\ \hline \end{tabular} \end{center}}

\allowdisplaybreaks
\begin{document}
\begin{center}

 {\Large \bf Unimodality of the Rank on Strongly Unimodal Sequences}
\end{center}

\begin{center}
  {Wenston J.T. Zang} \vskip 2mm

   School of Mathematics and Statistics, \& MOE Key Laboratory for Complexity Science in Aerospace, \& Xi'an-Budapest Joint Research Center for Combinatorics, Northwestern Polytechnical University, Xi'an 710072, P.R. China\\[6pt]

   \vskip 2mm

    zang@nwpu.edu.cn
\end{center}

\vskip 6mm \noindent {\bf Abstract.} Let $\{a_i\}_{i=1}^\ell$ be a strongly unimodal positive integer sequence with peak position $k$. The rank of such sequence is defined to be $\ell-2k+1$. Let $u(m,n)$ denote the number of sequences $\{a_i\}_{i=1}^\ell$ with rank $m$ and $\sum_{i=1}^{\ell} a_i=n$. Bringmann, Jennings-Shaffer, Mahlburg and Rhoades conjectured that  $\{u(m,n)\}_m$ is strongly log-concave for any fixed $n$. Motivated by this conjecture, in this paper we prove the strongly unimodality of $\{u(m,n)\}_m$, that is $u(m,n)>u(m+1,n)$ for $m\ge 0$ and $n\ge \max\{6,{m+2\choose 2}\}$. This result gives supportive evidence for the above conjecture. Moreover, we  find a combinatorial interpretation of $u(m,n)$, which leads to a new combinatorial interpretation of $\ospt(n)$. Furthermore, using this new combinatorial interpretation, a lower bound  and an asymptotic formula on $\ospt(n)$ will be presented.



\section{Introduction}

Let $(a_1,a_2,\ldots,a_\ell)$ be a strongly unimodal positive integer sequence of weight $n$, that is
\begin{equation}\label{ine-str-uni}
  1\le a_1<\cdots <a_{k-1}<a_k>a_{k+1}>\cdots>a_\ell\ge 1,
\end{equation}
for some $\ell\ge k\ge 1$ and $\sum_{i=1}^{\ell}a_i=n.$
Let $u(n)$ denote the number of such sequences. Andrews \cite{Andrews-2013} (where he use the word ``strictly convex composition" instead of ``strongly unimodal sequence") built the relation between the generating function of $u(n)$ and certain mock-theta functions. Bryson, Ono, Pitman and Rhoades \cite{Bryson-Ono-Pitman-Rhoades-2012} determined the parity of $u(n)$ on certain arithmetic progressions.  Rhoades \cite{Rhoades-2014} found the asymptotic formula on $u(n)$. Chen and Zhou \cite{Chen-Zhou-2016} related the generating function of $u(n)$ to a mixed mock modular form of weight $1/2$.

The rank of a strongly unimodal sequence have been extensively studied. Let $(a_1,\ldots,a_\ell)$ satisfies \eqref{ine-str-uni}, its rank is defined to be $\ell-2k+1$,
that is the number of integers after the maximum part minus the number of integers precedes the maximum part. Let $u(m,n)$ denote the number of strongly unimodal sequences of weight $n$ and rank $m$. Bryson, Ono, Pitman and Rhoades \cite{Bryson-Ono-Pitman-Rhoades-2012} gave the generating function of $u(m,n)$ as follows
\begin{equation}\label{equ-gen-umn}
  U(z;q)=\sum_{n=0}^{\infty}\sum_{m=-\infty}^{\infty}u(m,n)z^mq^n=
  \sum_{r=0}^{\infty}(-zq;q)_r(-q/z;q)_rq^{r+1}.
\end{equation}
Here we use the standard $q$-series notation
\[(a;q)_n=\prod_{i=1}^{n}(1-aq^{i-1})\quad\text{and}\quad
(a;q)_\infty=\prod_{i=1}^{\infty}(1-aq^{i-1}).\]
Moreover, they proved that $U(-1,q)$ is a quantum modular forms and $U(\pm i;q)$ is the third mock theta function.  There has also been related work connecting $U(z;q)$ with mock and quantum modular or Jacobi forms (see  \cite{Bringmann-Folsom-2016, Bringmann-Folsom-Rhoades-2015, Bringmann-Mahlburg-2014-jcta, Folsom-Ono-Rhoades-2013, Kim-Lim-Lovejoy-2016} for example).

The combinatorial properties of $u(m,n)$ was studied by
Bringmann, Jennings-Shaffer, Mahlburg and Rhoades \cite{Bringmann-Jennings-Mahlburg-Rhoades-2019}.
They gave a generating function of $u(m,n)$ for any fixed $m$. They also noticed that $\{u(m,n)\}_{m=-\infty}^{+\infty}$ appears to be unimodal. It is well known that for positive sequence $\{a_i\}$, the log-concavity of $\{a_i\}$ implies the unimodality of $\{a_i\}$. After checking by MAPLE for $n\le 500$, they raised the following conjecture on the strongly log-concavity of $u(m,n)$.

\begin{conj}[\cite{Bringmann-Jennings-Mahlburg-Rhoades-2019}]\label{conj-BJMR}
For $n\ge \max\{7,\frac{|m|(|m|+1)}{2}+1\}$,
  \[u(m,n)^2>u(m+1,n)u(m-1,n).\]
\end{conj}

%

The  asymptotic formulas on $u(m,n)$ was also given by  Bringmann, Jennings-Shaffer, Mahlburg and Rhoades \cite{Bringmann-Jennings-Mahlburg-Rhoades-2019} as follows.

\begin{thm}[\cite{Bringmann-Jennings-Mahlburg-Rhoades-2019}]
  For fixed $m\ge 0$, as $n\rightarrow\infty$, we have
  \begin{align}\label{equ-sim-umn}
    u(m,n) &\sim \frac{1}{16\sqrt{3}n}e^{\pi\sqrt{\frac{2n}{3}}}, \\
    u(m,n)-u(m+1,n) & \sim \frac{\pi(2m+1)}{96\sqrt{2}n^{\frac{3}{2}}}e^{\pi\sqrt{\frac{2n}{3}}}, \label{eq-asy-umn}\\
    u(m,n)^2-u(m+1,n)u(m-1,n) &\sim \frac{\pi}{786\sqrt{6}n^{\frac{5}{2}}}e^{2\pi\sqrt{\frac{2n}{3}}}.\label{eq-asy-umn-logc}
  \end{align}
\end{thm}

It should be noted that \eqref{eq-asy-umn-logc} implies Conjecture \ref{conj-BJMR} holds for  fixed $m$ and sufficiently  large $n$.

The main purpose of this paper is to establish the unimodality and strongly unimodality of $u(m,n)$.

\begin{thm}\label{thm-main}
  For any $m\ge 0$,
  \begin{equation}\label{equ-thm-main}
    u(m,n)\ge u(m+1,n).
  \end{equation}
Moreover, for any $m\ge 0$ and $n\ge \max\{6,{m+2\choose 2}\}$,
 \begin{equation}\label{equ-thm-main-1}
    u(m,n)> u(m+1,n).
  \end{equation}
\end{thm}

By \eqref{equ-gen-umn}, it is clear that $u(m,n)=u(-m,n)$. Thus Theorem \ref{thm-main} implies that for any $n\ge 0$,
\[\cdots\le u(-2,n)\le u(-1,n)\le u(0,n)\ge u(1,n)\ge u(2,n)\ge \cdots.\]
Thus $\{u(m,n)\}_{m=-\infty}^{+\infty}$ is unimodal. Moreover, for fixed $m\ge 0$ and
 $n\ge \max\{6,{m+2\choose 2}\}$,
\[u(-m-1,n)<\cdots < u(-1,n)< u(0,n)> u(1,n)> \cdots>u(m+1,n),\]
which implies the strongly unimodality of $\{u(k,n)\}_{k=-m-1}^{m+1}$.

Bringmann, Jennings-Shaffer, Mahlburg and Rhoades \cite{Bringmann-Jennings-Mahlburg-Rhoades-2019} asked if there are other combinatorial interpretations of $u(m,n)$.
In this paper, we find that $u(m,n)$ can be combinatorially interpreted in terms of rank and rank-set on integer partitions.  We recall that $\lambda=(\lambda_1,\ldots,\lambda_\ell)$ is an integer partition of $n$ if $\sum_{i=1}^{\ell}\lambda_i=n$ and
$\lambda_1\ge \lambda_2\ge\cdots\ge\lambda_\ell\ge 1.$
For the rest of this paper, we use $\ell(\lambda)$ to denote the length of a partition $\lambda$.  Denote $|\lambda|$ to be the weight of $\lambda,$ namely $|\lambda|=\sum_{i=1}^{\ell(\lambda)}\lambda_i.$
We use  $s(\lambda)$ to denote the smallest part of $\lambda$, that is
$s(\lambda)=\lambda_{\ell(\lambda)}.$ Moreover, let $\lambda'=(\lambda_1',\ldots,\lambda'_{\lambda_1})$ be the conjugate partition of $\lambda$. That is for $1\le i\le \lambda_1$, define
\[\lambda'_i=\#\{j\colon 1\le j\le \ell(\lambda),\  \lambda_j\ge i\}.\]
It is clear that $\lambda'_1=\ell(\lambda)$, $\ell(\lambda')=\lambda_1$ and $(\lambda')'=\lambda$. For more information of conjugate, we refer to  \cite[Definition 1.8]{Andrews-1976}.
Furthermore we use the convention that $\lambda_k=0$ for any $k>\ell(\lambda)$. For example, let $\lambda=(5,5,3,3,2,1)$, then $\ell(\lambda)=6$, $|\lambda|=19$, $s(\lambda)=1$, $\lambda'=(6,5,4,2,2)$ and $\lambda_k=0$ for all $k\ge 7.$

We next recall the definition of rank and rank-set of a partition. Dyson \cite{Dyson-1944} introduced the rank of a partition $\lambda$ as the largest part of the partition minus the number of parts, namely,
    ${\rm rank}(\lambda)=\lambda_1-\ell(\lambda).$
 The rank-set of $\lambda$ was introduced by Dyson \cite{Dyson-1989} as an infinite sequence
$(-\lambda_1,1-\lambda_2,2-\lambda_3,\ldots,j-\lambda_{j+1},
\ldots,\ell(\lambda)-1-\lambda_{\ell(\lambda)},\ell(\lambda),\ell(\lambda)+1,\ldots).$ For instance, the rank of $(5,5,3,2,2,1)$ is $5-6=-1$ and the rank-set of $(5,5,3,2,2,1)$ equals $(-5,-4,-1,1,2,4,6,7,\ldots)$.

Now we may interpret $u(m,n)$   in terms of rank and rank-set of partitions.

\begin{thm}\label{thm-com-int-umn}
  For $m\ge 0$, $u(m,n)$ is the number of partitions $\lambda$ of $n$ such that the following three restrictions hold:
  \begin{itemize}
    \item[(1)]
    $${\rm rank}(\lambda)\le\begin{cases}
                              -m-1, & \mbox{if } m\ge 1 \\
                              0, & \mbox{if } m=0
                            \end{cases};$$
    \item[(2)] $m-1$ is in the rank-set of $\lambda$;
    \item[(3)] When $m\ge 2$, $\lambda_1>\lambda_2>\cdots>\lambda_m$.
  \end{itemize}
\end{thm}

From Theorem \ref{thm-com-int-umn}, we may deduce another combinatorial interpretation of $u(m,n)$. In fact,  for any partition $\lambda$ counted by $u(m,n)$ as stated in Theorem \ref{thm-com-int-umn}, define
$$\mu=\begin{cases}
        (\lambda_1-m+1,\lambda_2-m+2,\ldots,\lambda_{m-1}-1,
\lambda_m,\ldots,\lambda_{\ell(\lambda)}), & \mbox{if } m\ge 1;\\
        \lambda, & \mbox{if }m=0.
      \end{cases}$$
 It is easy to check that the following corollary is equivalent to Theorem \ref{thm-com-int-umn}.

\begin{core}\label{cor-com-in}
  For $m\ge 0$, $u(m,n)$ equals the number of partitions $\mu$ of $n-{m\choose 2}$ satisfies the following two restrictions:
  \begin{itemize}
    \item[(1)] ${\rm rank}(\mu)\le -2m$;
    \item[(2)] $m-1$ is in the rank-set of $\mu$.
  \end{itemize}
\end{core}

It should be noted that Corollary \ref{cor-com-in} seems simpler than Theorem \ref{thm-com-int-umn}. But our main purpose is to give a proof of Theorem \ref{thm-main}. It is more convenient for us to adopt  Theorem \ref{thm-com-int-umn} than Corollary \ref{cor-com-in}.


Our next result is related to the ospt function.
Recall that the ospt function was introduced by Andrews, Chan and Kim \cite{Andrews-Chan-Kim} as follows:
\begin{equation}\label{equ-def-ospt}
  {\rm ospt}(n)=\sum_{|\lambda|=n \atop {\rm crank}(\lambda)\ge 0}{\rm crank}(\lambda)-\sum_{|\lambda|=n \atop {\rm rank}(\lambda)\ge 0}{\rm rank}(\lambda),
\end{equation}
where  the crank of a partition was  defined by Andrews and Garvan \cite{Andrews-Garvan-1988} as the largest part if the partition contains no ones, and otherwise as the number of parts larger than the number of ones minus the number of ones.
Andrews, Chan and Kim \cite{Andrews-Chan-Kim} gave a combinatorial interpretation on ${\rm ospt}(n)$ in terms of even strings and odd strings, and hence they showed that
${\rm ospt}(n)>0.$
Bringmann and Mahlburg \cite{Bringmann-Mahlburg-2014} found the following asymptotic result on ${\rm ospt}(n)$.

\begin{thm}[\cite{Bringmann-Mahlburg-2014}]\label{thm-Bringmann-ospt-asy}
As $n\rightarrow\infty$,
\begin{equation}\label{Bringmann-ospt-asy}
 {\rm ospt}(n)\sim \frac{p(n)}{4}.
\end{equation}
\end{thm}

Mao \cite{Mao-2018} improved Theorem \ref{thm-Bringmann-ospt-asy} as follows.

\begin{thm}[\cite{Mao-2018}]\label{thm-mao-ospt}
  As $n\rightarrow \infty$,
  \[\ospt(n)=\frac{p(n)}{4}(1+O(n^{-\frac{1}{4}})).\]
\end{thm}

Let $N(m,n)$ denote the number of partitions of $n$ with rank $m$, and $M(m,n)$ denote the number of partitions of $n$ with crank $m$. Chan and Mao \cite{Chan-Mao-2014} proved inequalities related to ${\rm ospt}(n)$ as given below.

\begin{thm}[\cite{Chan-Mao-2014}]
We have
\begin{align}
\ospt(n)&>\frac{p(n)}{4}+\frac{N(0,n)}{2}-\frac{M(0,n)}{4},&\text{for }n\geq 8,\label{chan-mao-lowerbound}\\[3pt]
\ospt(n)&<\frac{p(n)}{4}+\frac{N(0,n)}{2}-\frac{M(0,n)}{4}+\frac{N(1,n)}{2},&\text{for }n\geq 7,\label{chan-mao-upbound}\\[3pt]
\ospt(n)&<\frac{p(n)}{2},&\text{for }n\geq 3. \label{ospt-up-n/2}
\end{align}
\end{thm}

Chen, Ji and the author \cite{Chen-Ji-Zang-2017} gave another combinatorial interpretation of $\ospt(n)$ in terms of certain reordering $\tau_n$ on the set of partitions of $n$. They also improved \eqref{ospt-up-n/2} as follows:

\begin{thm}[\cite{Chen-Ji-Zang-2017}]\label{thm-CJZ-B}
For $n\ge 2$,
  \begin{equation}\label{bringmann-ospt-pn2}
    \ospt(n)\le \frac{p(n)}{2}-\frac{M(0,n)}{2}.
  \end{equation}
\end{thm}

 Bringmann, Jennings-Shaffer, Mahlburg and Rhoades \cite{Bringmann-Jennings-Mahlburg-Rhoades-2019} found that
  \begin{equation}\label{equ-cor-ospt-u0n}
    u(0,n)={\rm ospt}(n).
  \end{equation}
  Using this connection, they reproved Theorem \ref{thm-CJZ-B}.

In this paper, we  give a combinatorial interpretation on $\ospt(n)$, which is a direct consequence of Theorem \ref{thm-com-int-umn} and \eqref{equ-cor-ospt-u0n}.

\begin{core}\label{core-ospt-new}
  For $n\ge 0$, ${\rm ospt}(n)$ counts the number of partitions $\lambda$ of $n$ such that ${\rm rank}(\lambda)\le 0$ and $-1$ is in the rank-set of $\lambda$.
\end{core}

We note that  Corollary \ref{core-ospt-new} implies Theorem \ref{thm-CJZ-B}, which will be demonstrated in Section \ref{2}. Moreover, using Corollary \ref{core-ospt-new}, we find the following lower bound of ${\rm ospt}(n)$.

\begin{thm}\label{cor-ospt-lower}
For $n\ge 6$,
\begin{equation}\label{ospt-lower}
  {\rm ospt}(n)> \frac{p(n)}{4}+\frac{3N(0,n)}{8}.
\end{equation}
\end{thm}

It is hard to tell \eqref{ospt-lower} and \eqref{chan-mao-lowerbound} which one is stronger. By \cite[Theorem 1.4]{Mao-2014},
\begin{equation}\label{Mao-asy}
  N(0,n)\sim M(0,n)\sim \frac{\pi}{4\sqrt{6n}}p(n)\sim \frac{\pi}{48\sqrt{2}n^{\frac{3}{2}}}e^{\pi\sqrt{\frac{2n}{3}}}.
\end{equation}
We see that \eqref{ospt-lower} improves \eqref{chan-mao-lowerbound} asymptotically.

Moreover, we have the following asymptotical result on $\ospt(n)$, which improves both Theorem \ref{thm-Bringmann-ospt-asy} and Theorem \ref{thm-mao-ospt}.
\begin{thm}\label{thm-asy-ospt-pn4}
  As $n\rightarrow\infty$,
  \begin{equation}\label{ospt-ays}
   \ospt(n)-\frac{p(n)}{4}\sim \frac{N(0,n)}{2}.
  \end{equation}
\end{thm}

This paper is organized as follows. In Section \ref{2}, we prove Theorem \ref{thm-com-int-umn} in terms of the Durfee symbol and the $m$-Durfee rectangle symbol, as introduced in \cite{Andrews-2007} and \cite{Chen-Ji-Zang-2015} respectively.   Section \ref{3} is devoted to prove Theorem \ref{thm-main} by considering either $m\ge 1$ or $m=0$.   In Section \ref{5}, we give a proof of Theorem \ref{cor-ospt-lower} and Theorem \ref{thm-asy-ospt-pn4} by building a connection between $\ospt(n)$ and $u(0,n)-u(1,n)$.  An open problem will be discussed in Section \ref{6}.

\section{The combinatorial interpretation of $u(m,n)$}\label{2}

In this section, we give a proof of Theorem \ref{thm-com-int-umn}. To this end, we first recall the definition of $m$-Durfee rectangle symbol, which was introduced by Chen, Ji and the author \cite{Chen-Ji-Zang-2015}. We then restate Theorem \ref{thm-com-int-umn} in terms of $m$-Durfee rectangle symbol, which is Theorem \ref{thm-res-com-int-umn}.  Next we give a proof of Theorem \ref{thm-res-com-int-umn}, in this way we verify Theorem \ref{thm-com-int-umn}. Finally, we show that Theorem \ref{thm-CJZ-B} can be deduced from  Corollary \ref{core-ospt-new}.

%

We first    recall the definition of  $m$-Durfee rectangle symbol. Let $\lambda$ be a partition. The Ferrers diagram of $\lambda$ is obtained by drawing a left-justified array of $n$ dots with $\lambda_i$ dots in the $i$th row.  The $m$-Durfee rectangle of $\lambda$ is defined to be the largest $(m+j)\times j$ rectangle contained in the Ferrers diagram of $\lambda$, see Gordon and Houten \cite{Gordon-Houten-1968}. In other words, the integer $j$ in the $m$-Durfee rectangle $(m+j)\times j$ can be defined as follows:
\begin{equation}\label{eq-def-m-dur}
  j=\begin{cases}
      \max\{k\colon 1\le k\le \ell(\lambda),\  \lambda_{k+m}\ge k\}, & \mbox{if } \ell(\lambda)\ge m+1; \\
      0, & \mbox{if } \ell(\lambda)\le m.
    \end{cases}
\end{equation} 
 The $m$-Durfee rectangle symbol of $\lambda$ is defined as
\begin{equation}\label{m-dur}
(\alpha,\beta)_{(m+j)\times j}=\left(\begin{array}{c}
           \alpha \\
           \beta
         \end{array}\right)_{(m+j)\times j}=\left(\begin{array}{cccc}
\alpha_1,&\alpha_2,&\ldots,&\alpha_s\\[3pt]
\beta_1,&\beta_2,&\ldots,&\beta_t
\end{array} \right)_{(m+j)\times j},
\end{equation}
where $(m+j)\times j$ is the $m$-Durfee rectangle of $\lambda$, the partition $\alpha$ is the columns to the right of the $m$-Durfee rectangle, and $\beta$ consists of rows below the $m$-Durfee rectangle. That is
\begin{equation}\label{equ-def-alpha}
  \alpha=(\lambda_1-j,\lambda_2-j,\ldots,\lambda_{m+j}-j)'
\end{equation}
and
\begin{equation}\label{equ-def-beta}
  \beta=(\lambda_{m+j+1},\ldots,\lambda_{\ell(\lambda)}).
\end{equation}
It should be noted that when $\ell(\lambda)\le m$, by \eqref{eq-def-m-dur}, we have $j=0$. Thus $\beta=\emptyset$ and $\alpha=\lambda'$. To be specific, the $m$-Durfee rectangle symbol of $\lambda$ takes the form
\begin{equation}\label{eq-def-m-dur-0}
\lambda=(\alpha,\beta)_{m\times 0}=(\lambda',\emptyset)_{m\times 0}.
\end{equation}

It is clear that there is a one-to-one correspondence between the set of partitions and the set of $m$-Durfee rectangle symbols. For the rest of this paper, we do not distinguish a partition $\lambda$ and its $m$-Durfee rectangle symbol $(\alpha,\beta)_{(m+j)\times j}$. Write
\begin{equation*}
\lambda=(\alpha,\beta)_{(m+j)\times j}=\left(\begin{array}{cccc}
\alpha\\[3pt]
\beta
\end{array} \right)_{(m+j)\times j},
\end{equation*}
we also consider the $m$-Durfee rectangle symbol as a partition.

For example, as shown in Figure \ref{fig-m-Durfee}, the $1$-Durfee rectangle symbol of $(11,10,10,8,6,\break 5,4,4,3,3,1)$ is
\begin{equation*}
  \left(\begin{array}{cccccc}
          5, & 4, & 4, & 3, & 3, & 1 \\
          4, & 4, & 3, & 3, & 1 &
        \end{array}\right)_{6\times 5}.
\end{equation*}
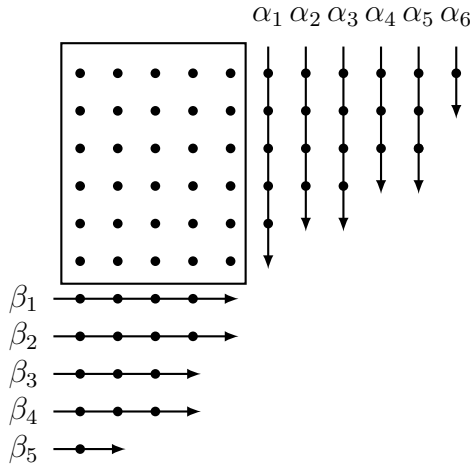
\begin{figure}[t]
  \centering
  \usetikzlibrary{arrows}
\begin{tikzpicture}
[place/.style={thick,fill=black!100,circle,inner sep=0pt,minimum size=1mm,draw=black!100}, sample/.style={thick,fill=black!100,circle,inner sep=0pt,minimum size=1mm,draw=black!100}]
\node [place] at (-3.5,2) {};
\node [place] (v7) at (-2,2) {};
\node [place] (v11) at (-1.5,2) {};
\node [place] (v16) at (-1,2) {};
\node [place] (v20) at (-0.5,2) {};
\node [place] (v24) at (0,2) {};
\node [place] at (-4,2) {};
\node [place] at (-5,2) {};
\node [place] at (-4.5,2) {};
\node [place] at (-3,2) {};
\node [place] (v5) at (-2.5,2) {};
\node [place] at (-5,1.5) {};
\node [place] at (-4.5,1.5) {};
\node [place] at (-4,1.5) {};
\node [place] at (-3.5,1.5) {};
\node [place] at (-3,1.5) {};
\node [place] (v4) at (-2.5,1.5) {};
\node [place] (v8) at (-2,1.5) {};
\node [place] (v12) at (-1.5,1.5) {};
\node [place] (v17) at (-1,1.5) {};
\node [place] (v21) at (-0.5,1.5) {};
\node [place] at (-5,1) {};
\node [place] at (-4.5,1) {};
\node [place] at (-4,1) {};
\node [place] at (-3.5,1) {};
\node [place] at (-3,0.5) {};
\node [place] at (-3,1) {};
\node [place] (v3) at (-2.5,1) {};
\node [place] (v9) at (-2,1) {};
\node [place] (v13) at (-1.5,1) {};
\node [place] (v18) at (-1,1) {};
\node [place] at (-0.5,1) {};
\node [place] at (-5,0.5) {};
\node [place] at (-4.5,0.5) {};
\node [place] at (-4,0.5) {};
\node [place] at (-3.5,0.5) {};
\node [place] (v2) at (-2.5,0.5) {};
\node [place] (v10) at (-2,0.5) {};
\node [place] (v14) at (-1.5,0.5) {};
\node [place] (v22) at (-0.5,1) {};
\node [place] at (-5,0) {};
\node [place] at (-4.5,0) {};
\node [place] at (-4,0) {};
\node [place] at (-3.5,0) {};
\node [place] at (-3,0) {};
\node [place] (v1) at (-2.5,0) {};
\node [place] at (-5,-0.5) {};
\node [place] at (-4.5,-0.5) {};
\node [place] at (-4,-0.5) {};
\node [place] at (-3.5,-0.5) {};
\node [place] at (-3,-0.5) {};
\node [place] (v26) at (-5,-1) {};
\node [place] (v27) at (-4.5,-1) {};
\node [place] (v28) at (-4,-1) {};
\node [place] (v29) at (-3.5,-1) {};
\node [place] (v30) at (-5,-1.5) {};
\node [place] (v31) at (-4.5,-1.5) {};
\node [place] (v32) at (-4,-1.5) {};
\node [place] (v33) at (-3.5,-1.5) {};
\node [place] (v34) at (-5,-2) {};
\node [place] (v35) at (-4.5,-2) {};
\node [place] (v36) at (-4,-2) {};
\node [place] (v37) at (-5,-2.5) {};
\node [place] (v38) at (-4.5,-2.5) {};
\node [place] (v39) at (-4,-2.5) {};
\node [place] (v40) at (-5,-3) {};
\draw  [thick](-5.25,2.4) rectangle (-2.8,-0.8);
\node (v6) at (-2.5,-0.75) {};
\node (v15) at (-1.5,-0.25) {};
\node (v19) at (-1,0.25) {};
\node (v23) at (-0.5,0.25) {};
\node (v25) at (0,1.25) {};
\node at (-2.5,2.75) {$\alpha_1$};
\node at (-2,2.75) {$\alpha_2$};
\node at (-1.5,2.75) {$\alpha_3$};
\node at (-1,2.75) {$\alpha_4$};
\node at (-0.5,2.75) {$\alpha_5$};
\node at (0,2.75) {$\alpha_6$};
\node at (-5.75,-1) {$\beta_1$};
\node at (-5.75,-1.5) {$\beta_2$};
\node at (-5.75,-2) {$\beta_3$};
\node at (-5.75,-2.5) {$\beta_4$};
\node at (-5.75,-3) {$\beta_5$};
\node (v49) at (-2.5,2.5) {};
\node (v50) at (-2,2.5) {};
\node (v51) at (-1.5,2.5) {};
\node (v52) at (-1,2.5) {};
\node (v53) at (-0.5,2.5) {};
\node (v54) at (0,2.5) {};
\node (v43) at (-5.5,-1) {};
\node (v44) at (-2.75,-1) {};
\node (v45) at (-5.5,-1.5) {};
\node (v46) at (-5.5,-2) {};
\node (v47) at (-5.5,-2.5) {};
\node (v48) at (-5.5,-3) {};
\draw[thick] [-latex](v43) edge (v44);
\node (v41) at (-2.75,-1.5) {};
\draw [thick] [-latex] (v45) edge (v41);
\node (v42) at (-3.25,-2) {};
\node (v55) at (-3.25,-2.5) {};
\node (v56) at (-4.25,-3) {};
\draw [thick] [-latex] (v46) edge (v42);
\draw [thick] [-latex] (v47) edge (v55);
\draw [thick] [-latex] (v48) edge (v56);
\node (v57) at (-2,-0.25) {};
\draw [thick] [-latex] (v49) edge (v6);
\draw [thick] [-latex] (v50) edge (v57);
\draw [thick] [-latex] (v51) edge (v15);
\draw [thick] [-latex] (v52) edge (v19);
\draw [thick] [-latex] (v53) edge (v23);
\draw [thick] [-latex] (v54) edge (v25);
\end{tikzpicture}
  \caption{An illustration of $1$-Durfee rectangle symbol of $(11,10,10,8,6,5,4,4,3,3,1)$.}\label{fig-m-Durfee}
\end{figure}

\begin{figure}[h]
  \centering
  \usetikzlibrary{arrows}
\begin{tikzpicture}
[place/.style={thick,fill=black!100,circle,inner sep=0pt,minimum size=1mm,draw=black!100}, sample/.style={thick,fill=black!100,circle,inner sep=0pt,minimum size=1mm,draw=black!100}]

\node [place] at (-3.5,2) {};
\node [place] (v7) at (-2,2) {};
\node [place] (v11) at (-1.5,2) {};
\node [place] (v16) at (-1,2) {};
\node [place] (v20) at (-0.5,2) {};
\node [place] (v24) at (0,2) {};
\node [place] at (-4,2) {};
\node [place] at (-5,2) {};
\node [place] at (-4.5,2) {};
\node [place] at (-3,2) {};
\node [place] (v5) at (-2.5,2) {};
\node [place] at (-5,1.5) {};
\node [place] at (-4.5,1.5) {};
\node [place] at (-4,1.5) {};
\node [place] at (-3.5,1.5) {};
\node [place] at (-3,1.5) {};
\node [place] (v4) at (-2.5,1.5) {};
\node [place] (v8) at (-2,1.5) {};
\node [place] (v12) at (-1.5,1.5) {};
\node [place] (v17) at (-1,1.5) {};
\node [place] (v21) at (-0.5,1.5) {};
\node [place] at (-5,1) {};
\node [place] at (-4.5,1) {};
\node [place] at (-4,1) {};
\node [place] at (-3.5,1) {};
\node [place] at (-3,0.5) {};
\node [place] at (-3,1) {};
\node [place] (v3) at (-2.5,1) {};
\node [place] (v9) at (-2,1) {};
\node [place] (v13) at (-1.5,1) {};
\node [place] (v18) at (-1,1) {};
\node [place] at (-0.5,1) {};
\node [place] at (-5,0.5) {};
\node [place] at (-4.5,0.5) {};
\node [place] at (-4,0.5) {};
\node [place] at (-3.5,0.5) {};
\node [place] (v2) at (-2.5,0.5) {};
\node [place] (v10) at (-2,0.5) {};
\node [place] (v14) at (-1.5,0.5) {};
\node [place] (v22) at (-0.5,1) {};
\node [place] at (-5,0) {};
\node [place] at (-4.5,0) {};
\node [place] at (-4,0) {};
\node [place] at (-3.5,0) {};
\node [place] at (-3,0) {};
\node [place] (v1) at (-2.5,0) {};
\node [place] at (-5,-0.5) {};
\node [place] at (-4.5,-0.5) {};
\node [place] at (-4,-0.5) {};
\node [place] at (-3.5,-0.5) {};
\node [place] at (-3,-0.5) {};
\node [place] (v26) at (-5,-1) {};
\node [place] (v27) at (-4.5,-1) {};
\node [place] (v28) at (-4,-1) {};
\node [place] (v29) at (-3.5,-1) {};
\node [place] (v30) at (-5,-1.5) {};
\node [place] (v31) at (-4.5,-1.5) {};
\node [place] (v32) at (-4,-1.5) {};
\node [place] (v33) at (-3.5,-1.5) {};
\node [place] (v34) at (-5,-2) {};
\node [place] (v35) at (-4.5,-2) {};
\node [place] (v36) at (-4,-2) {};
\node [place] (v37) at (-5,-2.5) {};
\node [place] (v38) at (-4.5,-2.5) {};
\node [place] (v39) at (-4,-2.5) {};
\node [place] (v40) at (-5,-3) {};
\draw  [thick](-5.25,2.4) rectangle (-2.8,-0.25);
\node (v6) at (-2.5,-0.5) {};
\node at (-2.5,2.75) {$\alpha_1$};
\node at (-2,2.75) {$\alpha_2$};
\node at (-1.5,2.75) {$\alpha_3$};
\node at (-1,2.75) {$\alpha_4$};
\node at (-0.5,2.75) {$\alpha_5$};
\node at (0,2.75) {$\alpha_6$};
\node at (-5.75,-0.5) {$\beta_1$};
\node at (-5.75,-1) {$\beta_2$};
\node at (-5.75,-1.5) {$\beta_3$};
\node at (-5.75,-2) {$\beta_4$};
\node at (-5.75,-2.5) {$\beta_5$};
\node at (-5.75,-3) {$\beta_6$};
\node (v49) at (-2.5,2.5) {};
\node (v50) at (-2,2.5) {};
\node (v51) at (-1.5,2.5) {};
\node (v52) at (-1,2.5) {};
\node (v53) at (-0.5,2.5) {};
\node (v54) at (0,2.5) {};
\node (v43) at (-5.5,-1) {};
\node (v44) at (-2.75,-1) {};
\draw [thick] [-latex] (v43) edge (v44);
\node (v45) at (-5.5,-1.5) {};
\node (v46) at (-5.5,-2) {};
\node (v47) at (-5.5,-2.5) {};
\node (v48) at (-5.5,-3) {};
\node (v41) at (-5.5,-0.5) {};
\draw [thick][-latex] (v41) edge (v6);
\draw [thick][-latex] (v49) edge (v6);
\node (v15) at (-2,-0.25) {};
\draw [thick][-latex] (v50) edge (v15);
\node (v19) at (-1.5,-0.25) {};
\node (v23) at (-1,0.25) {};
\node (v25) at (-0.5,0.25) {};
\node (v42) at (0,1.25) {};
\draw [thick][-latex] (v51) edge (v19);
\draw [thick][-latex] (v52) edge (v23);
\draw [thick][-latex] (v53) edge (v25);
\draw [thick][-latex] (v54) edge (v42);
\node (v55) at (-2.75,-1.5) {};
\node (v56) at (-3.25,-2) {};
\node (v57) at (-3.25,-2.5) {};
\node (v58) at (-4.25,-3) {};
\draw [thick][-latex] (v45) edge (v55);
\draw [thick][-latex] (v46) edge (v56);
\draw [thick][-latex] (v47) edge (v57);
\draw [thick][-latex] (v48) edge (v58);
\end{tikzpicture}
  \caption{An illustration of Durfee symbol of $(11,10,10,8,6,5,4,4,3,3,1)$.}\label{fig-Durfee}
\end{figure}
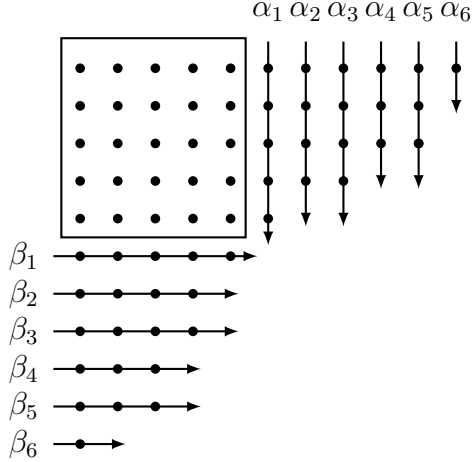

When $m=0$,  an $0$-Durfee rectangle is referred  to as a Durfee square, see \cite[P.28]{Andrews-1976}. The $0$-Durfee rectangle symbol is coincide with the Durfee symbol, see Andrews \cite{Andrews-2007}. In the notation of Andrews, a $d\times d$ square is simply denoted by $d$. For instance, as shown in Figure \ref{fig-Durfee}, the Durfee symbol of $(11,10,10,8,6,5,4,4,3,3,1)$ is
\[\left(\begin{array}{cccccc}
    5, & 4, & 4, & 3, & 3, & 1 \\
    5, & 4, & 4, & 3, & 3, & 1
  \end{array}\right)_5.\]

 Clearly, for any $m$-Durfee rectangle symbol \eqref{m-dur}, we see that $m+j\ge \alpha_1\ge\cdots\ge \alpha_s$, $j\ge \beta_1\ge\cdots\ge \beta_t$ and
\begin{equation}\label{equ-sum-m-drs}
  |\lambda|=|\alpha|+|\beta|+(m+j)j.
\end{equation}
Moreover, it is easy to see that
\begin{equation}\label{equ-rank-m-drs}
  {\rm rank}(\lambda)=\lambda_1-\ell(\lambda)=
  j+\ell(\alpha)-(m+j)-\ell(\beta)=\ell(\alpha)-\ell(\beta)-m.
\end{equation}

%

By \eqref{equ-rank-m-drs}, the following lemma is clear.

\begin{lem}\label{lem-equ-1}
  Let $\lambda$ be a partition of $n$ with its $m$-Durfee rectangle symbol $(\alpha,\beta)_{(m+j)\times j}$. Then the following two conditions are equivalent:
  \begin{itemize}
    \item[(1)]$${\rm rank}(\lambda)\le \begin{cases}
          -m-1, & \mbox{if } m\ge 1; \\
          0, & \mbox{if } m=0.
        \end{cases}$$
    \item[(2)]$$\ell(\alpha)-\ell(\beta)\le \begin{cases}
          -1, & \mbox{if } m\ge 1 ;\\
          0, & \mbox{if } m=0.
        \end{cases}$$
  \end{itemize}
\end{lem}

Evidently, Lemma \ref{lem-equ-1} restate the first restriction in Theorem \ref{thm-com-int-umn}  in terms of $m$-Durfee rectangle symbol. We next restate the second and the third restrictions in Theorem \ref{thm-com-int-umn} in terms of $m$-Durfee rectangle symbol as Lemma \ref{lem-equ-2} and Lemma \ref{lem-equ-3} respectively.

\begin{lem}\label{lem-equ-2}
  Let $\lambda$ be a partition of $n$ with its $m$-Durfee rectangle symbol $(\alpha,\beta)_{(m+j)\times j}$. Then the following two conditions are equivalent:
  \begin{itemize}
    \item[(1)] $m-1$ is in the rank-set of $\lambda$;
    \item[(2)] $\alpha_1\le m+j-1$.
  \end{itemize}
\end{lem}

\pf  On the one hand, if $\alpha_1\le m+j-1$, then by definition we see that $\lambda_{m+j}=j$, which implies $m+j-1-\lambda_{m+j}=m-1$ contains in the rank-set of $\lambda$.

On the other hand, if $m-1$ is in the rank-set of $\lambda$, then there exists $k$ such that $k-\lambda_{k+1}=m-1$. Hence either $\lambda_m=0$ or both $k\ge m$ and $\lambda_{k+1}=k-m+1$.

When $\lambda_m=0$, in other words $\ell(\lambda)<m$, by \eqref{eq-def-m-dur} we see that $j=0$ and from \eqref{eq-def-m-dur-0} we deduce that $\alpha_1=\ell(\lambda)\le m-1=m+j-1$.

When $k\ge m$ and $\lambda_{k+1}=k-m+1$, by \eqref{eq-def-m-dur} we see that the $m$-Durfee rectangle of $\lambda$ is $(k+1)\times (k-m+1)$. Thus $j=k-m+1$ and $\alpha_1=\#\{i\colon \lambda_i-(k-m+1)\ge 1\}$. From $\lambda_{k+1}=k-m+1$ we deduce that $\alpha_1\le k=m+j-1$. So in either case, we have $\alpha_1\le m+j-1$. This completes the proof.\qed

\begin{lem}\label{lem-equ-3}
Let $\lambda$ be a partition of $n$ with its $m$-Durfee rectangle symbol $(\alpha,\beta)_{(m+j)\times j}$. For any $1\le k\le m-1$, the following two conditions are equivalent:
  \begin{itemize}
    \item[(1)] $\lambda_k>\lambda_{k+1}$;
    \item[(2)] $k$ appears in $\alpha$.
  \end{itemize}
\end{lem}

\pf    On the one hand, if $\lambda_k>\lambda_{k+1},$  denote $\lambda_k-j$ by $s$, then  $s=\lambda_k-j>\lambda_{k+1}-j$. Thus from \eqref{equ-def-alpha}, we see that
 \begin{equation*}
   \alpha_s=\#\{i\colon 1\le i\le m+j,\  \lambda_i-j\ge s\}
 \end{equation*}
 and
 \[k+1>\#\{i\colon 1\le i\le m+j, \  \lambda_i-j\ge s\}\ge k.\]
Hence  $\alpha_s=k$.

 On the other hand, if there exists an integer $s$ such that $\alpha_s=k$, then
$\#\{t\colon \alpha_t\ge k\}\ge s>\#\{t\colon \alpha_t\ge k+1\}.$
 By \eqref{equ-def-alpha}, we see that
 $\lambda_k-j\ge s>\lambda_{k+1}-j,$
which implies $\lambda_k>\lambda_{k+1}$.
\qed

Given $m\ge 0$, let $U(m,n)$ denote the set of partitions $(\alpha,\beta)_{(m+j)\times j}$ of $n$ satisfies the following restrictions:
  \begin{itemize}
    \item[(1)] $$\ell(\alpha)-\ell(\beta)\le \begin{cases}
          -1, & \mbox{if } m\ge 1 \\
          0, & \mbox{if } m=0.
        \end{cases};$$
    \item[(2)] $\alpha_1\le m+j-1$;
    \item[(3)] $1,2,\ldots,m-1$ appears in $\alpha$.
  \end{itemize}
From Lemma \ref{lem-equ-1}, Lemma \ref{lem-equ-2} and Lemma \ref{lem-equ-3}, we see that Theorem \ref{thm-com-int-umn} is equivalent to the following theorem.

\begin{thm}\label{thm-res-com-int-umn}
  For $m\ge 0$, $\# U(m,n)=u(m,n).$
\end{thm}

\pf  It sufficient to show that for fixed $m\ge 0$, the generating function of $u(m,n)$ equals the generating function of $\# U(m,n)$.  We first calculate the  generating function of $u(m,n)$ as follows.

Let
\[{r\brack j}_q=\frac{(q;q)_r}{(q;q)_j(q;q)_{r-j}}\]
be the Gauss polynomials. It is well known that the ${r\brack j}_q$ is the generating function of partitions with at most $j$ parts, each part not exceeding $r-j$ (see \cite[Theorem 3.1]{Andrews-1976}). To be specific,
\begin{equation}\label{gen-gaus}
  \sum_{\lambda_1\le r-j\atop \ell(\lambda)\le j}q^{|\lambda|}={r\brack j}_q.
\end{equation}
We also recall the following result (see \cite[Theorem 3.3]{Andrews-1976}).
\begin{equation}\label{eq-gaus}
(-zq;q)_k=\sum_{j=0}^{k}{k\brack j}_qz^jq^{j+1\choose 2}.
\end{equation}
Using \eqref{eq-gaus}, we may transform \eqref{equ-gen-umn} as given below.
\begin{equation}\label{equ-gen-umn-1}
  \sum_{n=0}^{\infty}\sum_{m=-\infty}^{\infty}u(m,n)z^mq^n=
  \sum_{r=0}^{\infty}\sum_{j=0}^{r}\sum_{k=0}^{r}z^{j-k}{r\brack j}_q{r\brack k}_qq^{r+1+{j+1\choose 2}+{k+1\choose 2}}.
\end{equation}
Note that $m\ge 0$,
comparing the coefficient of $z^m$ on both sides of \eqref{equ-gen-umn-1}, we have
  \begin{equation}\label{equ-gen-umn-2}
 \sum_{n=0}^{\infty}u(m,n)q^n=
  \sum_{r=m}^{\infty}\sum_{k=0}^{r-m}{r\brack k+m}_q{r\brack k}_qq^{r+1+{k+m+1\choose 2}+{k+1\choose 2}}.
\end{equation}

We next calculate the generating function of $\#U(m,n)$. We first consider the case  $m\ge 1$.
For any $(\alpha,\beta)_{(m+j)\times j}\in U(m,n)$, by $\ell(\alpha)-\ell(\beta)\le -1$, we see that $\ell(\beta)\ge 1$. Thus $j\ge \beta_1\ge 1$. Set $\ell(\beta)=t$, clearly $\ell(\alpha)\le t-1$. Moreover,  $1,2,\ldots,m-1$ is a part of $\alpha$ implies that $\ell(\alpha)\ge m-1$, therefore $t\ge \ell(\alpha)+1\ge m$.

The generating function of $\alpha$ can be determined as follows.
Remove $1,2,\ldots,m-1$ exactly once from $\alpha$ to obtain a new partition $\bar{\alpha}$.
For example, let $m=4$ and $\alpha=(7,7,6,4,3,2,2,1,1,1)$, then $\bar{\alpha}=(7,7,6,4,2,1,1)$.
Clearly $\ell(\bar{\alpha})=\ell(\alpha)-m+1\le t-m,$ and $\bar{\alpha}_1\le \alpha_1\le m+j-1$. By \eqref{gen-gaus},
\begin{equation}\label{gen-gauss}
  \sum_{\bar{\alpha}_1\le j+m-1\atop \ell(\bar{\alpha})\le t-m} q^{|\bar{\alpha}|}={t+j-1\brack j+m-1}_q.
\end{equation}
Hence for fixed $t\ge m$ and $j\ge 1$, the generating function of $\alpha$ equals
\begin{equation}\label{gen-alpha}
\sum_{\alpha} q^{|{\alpha}|}=\sum_{\bar{\alpha}_1\le j+m-1\atop \ell(\bar{\alpha})\le t-m} q^{|\bar{\alpha}|}q^{m\choose 2}={t+j-1\brack j+m-1}_qq^{m\choose 2}.
\end{equation}

We next establish the generating function of $\beta$. Define $\bar{\beta}=(\beta_1-1,\ldots,\beta_{t}-1)$. Thus $\bar{\beta}$ is a partition with each part not exceeding $j-1$ and $\ell(\bar{\beta})\le t$. From \eqref{gen-gaus}, for fixed $t\ge m$ and $j\ge 1$, the generating function of $\beta$ can be expressed as follows:
\begin{equation}\label{gen-beta}
\sum_{\beta} q^{|{\beta}|}=\sum_{\bar{\beta}_1\le j-1\atop \ell(\bar{\beta})\le t} q^{|\bar{\beta}|}q^{t}={t+j-1\brack j-1}_qq^{t}.
\end{equation}

Now we may calculate the generating function of $\#U(m,n)$.  By \eqref{equ-sum-m-drs}, \eqref{gen-alpha} and \eqref{gen-beta},
\begin{align}
  \sum_{n=0}^{\infty}\# U(m,n)q^n&=\sum_{t=m}^{\infty}\sum_{j=1}^{\infty} q^{j(m+j)}\sum_{\alpha}q^{|\alpha|}
  \sum_{\beta}q^{|\beta|}\label{subalpha-1}\\
  &=\sum_{t=m}^{\infty}\sum_{j=1}^{\infty} q^{j(m+j)}{t+j-1\brack j+m-1}_qq^{m\choose 2}{t+j-1\brack j-1}_qq^{t}\nonumber\\
&\xlongequal[r=t+k]{k=j-1}\sum_{r=m}^{\infty}\sum_{k=0}^{r-m} q^{(k+1)(m+k+1)}{r\brack k+m}_qq^{m\choose 2}{r\brack k}_qq^{r-k}\nonumber\\
  &=\sum_{r=m}^{\infty}\sum_{k=0}^{r-m}{r\brack k+m}_q{r\brack k}_qq^{r+1+{k+m+1\choose 2}+{k+1\choose 2}},\label{gen-am}
\end{align}
 where the summation in \eqref{subalpha-1} ranges over all the $\alpha$ and $\beta$ satisfies $\ell(\alpha)\le t-1$, $\ell(\beta)=t$, $\alpha_1\le m+j-1$, $\beta_1\le j$ and $1,2,\ldots,m-1$ appears in $\alpha$. Clearly, \eqref{gen-am} equals \eqref{equ-gen-umn-2}.

We now consider the case $m=0$.  For any partition $(\alpha,\beta)_d\in U(0,n)$, by definition, $\ell(\alpha)-\ell(\beta)\le 0$, $\beta_1\le d$ and $\alpha_1\le d-1$. Denote $\ell(\beta)=t$, which implies $t\ge \ell(\alpha)\ge 0$.
We first calculate the generating function of $\alpha$. By \eqref{gen-gaus},  for fixed $t\ge 0$ and $d\ge 1$,
\begin{equation}\label{gen-alpha-0}
  \sum_{\alpha_1\le d-1\atop \ell(\alpha)\le t}q^{|\alpha|}={t+d-1\brack d-1}_q.
\end{equation}

On the other hand,  set $\bar{\beta}=(\beta_1-1,\ldots,\beta_{t}-1)$. Then $\bar{\beta}_1\le d-1$ and $\ell(\bar{\beta})\le t$. From \eqref{gen-gaus}, for fixed $t\ge 0$ and $d\ge 1$,
\begin{equation}\label{gen-beta-0}
 \sum_{\beta_1\le d\atop \ell(\beta)=t} q^{|\beta|}= \sum_{\bar{\beta}_1\le d-1\atop \ell(\bar{\beta})\le t}q^{|\bar{\beta}|}q^{t}={t+d-1\brack d-1}_qq^{t}.
\end{equation}

Using \eqref{equ-sum-m-drs}, \eqref{gen-alpha-0} and \eqref{gen-beta-0},  we have
\begin{align}
  \sum_{n=0}^{\infty}\# U(0,n)q^n&=
  \sum_{d=1}^{\infty}\sum_{t=0}^\infty q^{d^2}
  \sum_{\alpha}q^{|\alpha|}\sum_{\beta}q^{|\beta|} \label{0-alpha-beta}\\
  &=\sum_{d=1}^{\infty}\sum_{t=0}^{\infty} {t+d-1\brack d-1}_q^2q^{d^2+t}\nonumber\\
 &\xlongequal[r=t+k]{k=d-1} \sum_{r=0}^{\infty}\sum_{k=0}^{r}{r\brack k}_q^2q^{k^2+k+r+1}\label{gen-a0n},
\end{align}
where the summations in \eqref{0-alpha-beta} ranges over all $\alpha,\beta$ satisfies $\ell(\alpha)\le \ell(\beta)=t$, $\alpha_1\le d-1$ and $\beta_1\le d$. Clearly when $m=0$, we have \eqref{gen-a0n} equals \eqref{equ-gen-umn-2}. \qed

As stated in Introduction, Corollary \ref{core-ospt-new} is a direct consequence of Theorem \ref{thm-com-int-umn}. Now we conclude this section by showing that Theorem \ref{thm-CJZ-B} can be deduced from Corollary \ref{core-ospt-new}.

{\noindent \it Proof of Theorem \ref{thm-CJZ-B}.} Let $q(-1,n)$ denote the number of partitions $\lambda$ of $n$ such that
$-1$ is in the rank-set of $\lambda$. Dyson \cite{Dyson-1989} (see also Berkovich and Garvan \cite{Berkovich-Garvan-2002}) showed that
\begin{equation}\label{eq-q1n-mmn}
  q(-1,n)=\sum_{m=-\infty}^{-1}M(m,n).
\end{equation}
From  $M(m,n)=M(-m,n)$ due to Dyson \cite{Dyson-1989}, we deduce that
\begin{equation}\label{eq-pn-mmn}
  p(n)=\sum_{m=-\infty}^{+\infty}M(m,n)=M(0,n)+2\sum_{m=-\infty}^{-1}M(m,n).
\end{equation}
Combining \eqref{eq-q1n-mmn} and \eqref{eq-pn-mmn}, we have
\begin{equation}\label{eq-q1-pn}
  q(-1,n)=\frac{p(n)-M(0,n)}{2}.
\end{equation}
On the other hand, by Corollary \ref{core-ospt-new}, ${\rm ospt}(n)$ counts certain subset of the set counted by $q(-1,n)$. Therefore
\[\ospt(n)\le q(-1,n).\]
Thus by \eqref{eq-q1-pn}, we derive our desired result.\qed

%

\section{The unimodality of $u(m,n)$}\label{3}

This section is devoted to prove Theorem \ref{thm-main}. By Theorem \ref{thm-res-com-int-umn}, the set $U(m+1,n)$ is defined in terms of $(m+1)$-Durfee rectangle symbols. In order to prove Theorem \ref{thm-main}, we first redefine the set $U(m+1,n)$ as $V(m+1,n)$ by using $m$-Durfee rectangle symbols instead of $m+1$-Durfee rectangle symbols. We then  present a proof of Theorem \ref{thm-main} for $m\ge 1$ by building an injection $\Phi$ from the set $V(m+1,n)$ into the set $U(m,n)$. When $m=0$, the case is more complicated. We aim to build an injection $\Psi$ from $V(1,n)$ into $U(0,n)$. To this end,   we divide the set $V(1,n)$ into three disjoint subsets $V_1(1,n)$, $V_2(1,n)$ and $V_3(1,n)$. Also we divide the set $U(0,n)$ into five disjoint subsets $U_i(0,n)$ for $1\le i\le 5$. The injection $\Psi$ consists three injections $\psi_i$ from $V_i(1,n)$ into $U_i(0,n)$ for $1\le i\le 3$. In this way, we verify Theorem \ref{thm-main}.

For $m\ge 0$, let $V(m+1,n)$ denote the set of partitions $\lambda$ with its $m$-Durfee rectangle symbol $(\alpha,\beta)_{(m+j)\times j}$ of $n$ satisfies the following three restrictions:
  \begin{itemize}
    \item[(1)] $\ell(\alpha)-\ell(\beta)\le
          -2;$
    \item[(2)] $\beta_1=j$;
    \item[(3)] $1,2,\ldots,m$ appears in $\alpha$.
  \end{itemize}

The following corollary gives the cardinality of $V(m+1,n)$.

\begin{core}\label{cor-res-com-int-umn}
  For $m\ge 0$,  $\# V(m+1,n)=u(m+1,n).$
\end{core}

\pf To show $\# V(m+1,n)=u(m+1,n)$, it sufficient to prove that the three restrictions in the definition of $V(m+1,n)$ are equivalent to the three restrictions in Theorem \ref{thm-com-int-umn} for $u(m+1,n)$ respectively. In fact, the equivalences of  the first and third restrictions can be proved by using the same argument as in the proof of Lemma \ref{lem-equ-1} and Lemma \ref{lem-equ-3}. The equivalence of the second restriction is a direct consequence of \cite[Proposition 3.1]{Chen-Ji-Zang-2015}. \qed

We are now in a position to prove Theorem \ref{thm-main} for $m\ge 1$.

{\noindent \it Proof of Theorem \ref{thm-main} for $m\ge 1$.} We aim to build an injection $\Phi$ from $V(m+1,n)$ into $U(m,n)$. This yields $u(m,n)\ge u(m+1,n)$ as desired.

Clearly, $V(m+1,n)$ can be divided into the following two disjoint subsets:
\begin{itemize}
  \item[(1)] $V_1(m+1,n)$ is the set of partitions $(\alpha,\beta)_{(m+j)\times j}\in V(m+1,n)$ such that $\alpha_1<m+j$;
  \item[(2)] $V_2(m+1,n)$ is the set of partitions $(\alpha,\beta)_{(m+j)\times j}\in V(m+1,n)$ such that $\alpha_1=m+j$.
\end{itemize}
Obviously $V_1(m+1,n)\cup V_2(m+1,n)=V(m+1,n)$ and
$V_1(m+1,n)\cap V_2(m+1,n)=\emptyset$.

On the other hand, we define two disjoint subsets of $U(m,n)$ as given below:
\begin{itemize}
  \item[(1)] $U_1(m,n)$ is the set of partitions $(\gamma,\delta)_{(m+j')\times j'}\in U(m,n)$ satisfies the following three restrictions:
      \begin{itemize}
  \item[(1)] $\ell(\gamma)-\ell(\delta)\le -2$;
  \item[(2)] $\delta_1=j'$;
  \item[(3)] $\gamma$ contains $m$ as a part.
\end{itemize}
  \item[(2)] $U_2(m,n)$ is the set of partitions $(\gamma,\delta)_{(m+j')\times j'}\in U(m,n)$ such that $\delta_1<j'$.
\end{itemize}
Evidently, $U_1(m,n)\cap U_2(m,n)=\emptyset$. We next build an injection $\phi_1$ from $V_1(m+1,n)$ into $U_1(m,n)$ and an injection $\phi_2$ from $V_2(m+1,n)$ into $U_2(m,n)$. The injection $\Phi$ can be defined as follows:
\begin{equation}\label{equ-def-Phi}
  \Phi(\lambda)=\begin{cases}
                  \phi_1(\lambda), & \mbox{if } \lambda\in V_1(m+1,n); \\
                  \phi_2(\lambda), & \mbox{if } \lambda\in V_2(m+1,n).
                \end{cases}
\end{equation}
When $\lambda\in V_1(m+1,n)$, it is trivial to check that $V_1(m+1,n)$ coincides with $U_1(m,n)$, thus the map $\phi_1$ can set to be the identity map.

We proceed to describe the map $\phi_2$. Let
\begin{equation}
\lambda=\left(\begin{array}{cccc}
\alpha_1,&\alpha_2,&\ldots,&\alpha_s\\[3pt]
\beta_1, &\beta_2,&\ldots,&\beta_t
\end{array} \right)_{(m+j)\times j}
\end{equation}
be a partition in $V_2(m+1,n)$. By definition, $s-t\le -2$, $\alpha_1=m+j$ and $\beta_1=j$. For any $1\le k\le m$, since $k$ appears in $\alpha$, we may define $i_k$ to be the unique integer such that $\alpha_{i_k}=k$ and $\alpha_{i_k+1}=k-1$. If $j=0$, then by \eqref{eq-def-m-dur-0} we have $t=0$, contradicts to $s-t\le -2$. Thus $j\ge 1$ and $\alpha_1=j+m>m$, hence $i_m\ge 2$. Define
\begin{equation*}
\phi_2(\lambda)=\left(\begin{array}{c}
                      \gamma \\
                      \delta
                    \end{array}\right)_{(m+j')\times j'}=\left(\begin{array}{cccccccc}
\alpha_2,&\ldots,&\alpha_{i_m-1},&\alpha_{i_m}-1,&\alpha_{i_m+1},&\ldots,&\alpha_s\\[3pt]
\beta_2, &\beta_3, &\ldots, &\beta_t&&&
\end{array} \right)_{(m+j+1)\times (j+1)}.
\end{equation*}
By the definition of $i_m$, we see that
$\alpha_{i_m}-1=m-1=\alpha_{i_m+1},$
which implies $\gamma$ is a partition. Thus $\phi_2(\lambda)$ is well defined. We next show that $\phi_2(\lambda)\in U_2(m,n)$. By definition, it sufficient to show that  $|\phi_2(\lambda)|=|\lambda|$, $\gamma_1\le m+j'-1$, $\ell(\gamma)-\ell(\delta)\le -1$, $\delta_1<j'$ and $1,2,\ldots,m-1$ appears in $\gamma$.

First,
\begin{align}\label{equ-def-phi}
  |\phi_2(\lambda)| &=|\gamma|+|\delta|+(m+j+1)(j+1) \nonumber\\
   &=\left(\sum_{i=2}^{s} \alpha_i\right)-1+\left(\sum_{i=2}^{t}\beta_i\right)+j(m+j)+j+m+j+1 \nonumber\\
   &=\left(m+j+\sum_{i=2}^{s} \alpha_i\right)+\left(j+\sum_{i=2}^{t}\beta_i\right)+j(m+j)
   \nonumber\\
   &=|\alpha|+|\beta|+(m+j)j=|\lambda|.
\end{align}

Moreover, it is clear that $\ell(\delta)=t-1$.  When $m=1$, note that $i_1=s$, we see that $\gamma=(\alpha_2,\ldots,\alpha_{s-1})$, which implies $\ell(\gamma)=s-2$. When $m\ge 2$, we have $\ell(\gamma)=s-1$. Define $\delta_{m,1}=1$ for $m=1$ and $\delta_{m,1}=0$ for $m\ne 1$, from the above analysis we have
\begin{equation}\label{equ-lga-ldel}
  \ell(\gamma)-\ell(\delta)=(s-1-\delta_{m,1})-(t-1)=s-t-\delta_{m,1}\le -2-\delta_{m,1}.
\end{equation}
Thus we verify $\ell(\gamma)-\ell(\delta)\le -1$.

Furthermore, it is clear that $\gamma_1= \alpha_2\le m+j=m+j'-1$, $\gamma_{i_k-1}=\alpha_{i_k}=k$ for $1\le k\le m-1$ and $\delta_1=\beta_2\le \beta_1=j<j'$. In this way we show that $\phi_2(\lambda)\in U_2(m,n)$.

It remains to show $\phi_2$ is an injection. Let
$I(m,n)=\{\phi_2(\lambda)\colon \lambda\in V_2(m+1,n)\}$
be the image set of $\phi_2$, which has been proved to be a subset of $U_2(m,n)$. We will construct a map $\varphi$ from $I(m,n)$ to $V_2(m+1,n)$, such that for any $\lambda\in V_2(m+1,n)$, we have
\begin{equation}\label{equ-def-inv}
  \varphi(\phi_2(\lambda))=\lambda.
\end{equation}

For any $\mu=(\gamma,\delta)_{(m+j')\times j'}\in I(m,n)\subseteq U_2(m,n)$, by definition we have $\delta_1\le j'-1$. Moreover, for any $1\le k\le m-1$, define $t_k=\min\{p\colon \gamma_p=k\}.$ Since $1,2,\ldots,m-1$ appears in $\gamma$, we see such $t_k$ exists.   We also use the convention that $t_0=\ell(\gamma)+1$. By \eqref{equ-lga-ldel}, we see that $\ell(\gamma)-\ell(\delta)\le -2-\delta_{m,1}.$
Furthermore, from the construction of $\phi_2$, we see that $\gamma_{i_{m}-1}=\alpha_{i_m}-1=m-1$ and $\gamma_{i_m}=\alpha_{i_m+1}=m-1$. Thus the part $m-1$ appears at least twice in $\gamma$. Equivalently,   $\gamma_{t_{m-1}}=\gamma_{t_{m-1}+1}=m-1$. We may define $\varphi(\mu)$ as follows:
\begin{align*}
\varphi(\mu)&=\left(\begin{array}{c}
                         \alpha \\
                         \beta
                       \end{array}\right)_{(m+j)\times j}\\
&=\left(\begin{array}{cccccccc}
m+j'-1,&\gamma_1,&\ldots,&\gamma_{t_{m-1}-1},&\gamma_{t_{m-1}}+1,&\gamma_{t_{m-1}+1},&\ldots,&\gamma_{\ell(\gamma)}\\[3pt]
j'-1, &\delta_1, &\ldots, &\delta_{\ell(\delta)}
\end{array} \right)_{(m+j'-1)\times (j'-1)}.
\end{align*}
From the minimum of $t_{m-1}$, we see that $\gamma_{t_{m-1}-1}\ge m=\gamma_{t_{m-1}}+1$. By definition, $\gamma_1\le m+j'-1$ and $\delta_1\le j'-1$. Hence $\alpha$ and $\beta$ are partitions, which yields $\varphi(\mu)$ is well defined.

We next show that $\varphi(\mu)\in V_2(m+1,n)$. Evidently, $\beta_1=j'-1=j$ and $\alpha_1=m+j'-1=m+j$. Moreover,  $\ell(\alpha)-\ell(\beta)=\ell(\gamma)-\ell(\delta)+\delta_{1,m}\le -2.$ Furthermore, $\alpha_{t_{m-1}+1}=m$, $\alpha_{t_{k}+1}=k$ for $1\le k\le m-1$ and it can be checked that $|\varphi(\mu)|=|\mu|$. Thus $\varphi(\mu)\in V_2(m+1,n)$. From the construction of $\phi_2$ and $\varphi$, we see that \eqref{equ-def-inv} holds. Thus  $\phi_2$ is an injection, which implies the map $\Phi$ as defined in \eqref{equ-def-Phi} is an injection.

We conclude this proof by showing that $u(m,n)>u(m+1,n)$ holds when $m\ge 1$ and $n\ge {m+2\choose 2}$. It sufficient to find a partition $\mu\in U(0,n)$ and $\mu\not\in U_1(0,n)\cup U_2(0,n)$. Let $n-{m+2\choose 2}=r\ge 0$. Define
\[\mu=\left(\begin{array}{c}
                  \gamma \\
                  \delta
                \end{array}\right)_{(m+j')\times j'}=\left(\begin{array}{ccccc}
            m-1, & m-2,&,\ldots,&1 \\
            1^{m+r} & &&
          \end{array}\right)_{(m+1)\times 1}.\]
Clearly, $\ell(\gamma)-\ell(\delta)=m-1-m-r\le -1$, $\gamma_1=m-1<(m+j')-1$ and $1,2,\ldots,m-1$ appears in $\gamma$. Moreover,
$|\mu|={m\choose 2}+m+r+m+1={m+2\choose 2}+r=n.$
Therefore $\mu\in U(m,n)$. From $\delta_1=j'=1$, we see that $\mu\not\in U_2(m,n)$. Since $m$ is not a part of $\gamma$, it follows that $\mu\not\in U_1(m,n)$. This completes the proof.\qed

For example, let $m=2$ and
$$\lambda=\left(\begin{array}{ccccccccccc}
                  5, & 5, &3,& 2, & 2, & 1 \\
                  3, & 3, & 3, & 2, & 2,&2,&1,&1,&1
                \end{array}\right)_{5\times 3}\in V(3,51).$$
We see that $j=3$ and $\alpha_1=5=m+j$, thus  $\lambda\in V_2(3,51)$. We need to apply $\phi_2$ on $\lambda$. It is clear that $i_2=5$ and
$$\phi_2(\lambda)=\mu=\left(\begin{array}{ccccccccccc}
                  5, & 3,&2, & 1, & 1 \\
                  3, & 3, & 2, & 2,&2,&1,&1,&1
                \end{array}\right)_{6\times 4}\in U_2(2,51).$$
Applying $\varphi$ in $\mu$, we see that $t_1=4$, and $\varphi(\mu)=\lambda$.

We next consider the case $m=0$. To build an injection $\Psi$ from $V(1,n)$ into $U(0,n)$, we will partition $V(1,n)$ into three blocks $V_i(1,n)$ $(1\le i\le 3)$ and partition $U(0,n)$ into five blocks $U_i(0,n)$ $(1\le i\le 5)$. The injection $\Psi$ consists of three injections $\psi_i$ from $V_i(1,n)$ into $U_i(0,n)$ for $1\le i\le 3$.
For any $(\alpha,\beta)_d\in V(1,n)$, by definition $\alpha_1\le d$. The blocks $V_i(1,n) (1\le i\le 3)$ can be described as follows:
\begin{itemize}
  \item[(1)] $V_1(1,n)$ is the subset of $V(1,n)$ such that $\alpha_1\le d-1$;
  \item[(2)] $V_2(1,n)$ is the subset of   $V(1,n)$ such that $\alpha_1=d$ and $s(\beta)=1$;
  \item[(3)]$V_3(1,n)$ is the subset of   $V(1,n)$ such that $\alpha_1=d$ and $s(\beta)\ge 2$.
\end{itemize}
 It is easy to check that $\{V_1(1,n),V_2(1,n), V_3(1,n)\}$ is a set partition of $V(1,n)$.

Given $(\gamma,\delta)_{ d'}\in U(0,n)$, by definition $\delta_1\le d'$ and $\ell(\gamma)-\ell(\delta)\le 0$. The blocks $U_i(0,n)$ $(1\le i\le 5)$ can be defined as given below.
\begin{itemize}
  \item[(1)] $U_1(0,n)$ is the subset of $U(0,n)$ such that $\ell(\gamma)-\ell(\delta)\le -2$ and $\delta_1=d'$;
  \item[(2)] $U_2(0,n)$ is the subset of $U(0,n)$ such that $\ell(\gamma)-\ell(\delta)\le -1$ and $\delta_1\le d'-1$;
  \item[(3)] $U_3(0,n)$ is the subset of $U(0,n)$ such that $\ell(\gamma)-\ell(\delta)=0$ and $\delta_1\le d'-1$;
  \item[(4)] $U_4(0,n)$ is the subset of $U(0,n)$ such that
          $\ell(\gamma)-\ell(\delta)=-1$ and $\delta_1=d'$;
  \item[(5)] $U_5(0,n)$ is the subset of $U(0,n)$ such that $\ell(\gamma)-\ell(\delta)=0$ and $\delta_1= d'$.
\end{itemize}

It is trivial to check that $\{U_1(0,n),U_2(0,n),U_3(0,n),U_4(0,n),U_5(0,n)\}$ is a set partition of $U(0,n)$. We are now ready to present three injections $\psi_i$ from $V_i(1,n)$ into $U_i(0,n)$ ($1\le i\le 3$). It is clear that $V_1(1,n)$ is coincide with $U_1(0,n)$. Thus $\psi_1$ can be set to the identity map. The following lemma gives a bijection $\psi_2$ between $V_2(1,n)$ and $U_2(0,n)$.

\begin{lem}\label{lem-psi2}
  There is a bijection $\psi_2$ between $V_2(1,n)$ and $U_2(0,n)$.
\end{lem}

\pf Let
\begin{equation}
\lambda=\left(\begin{array}{c}
                        \alpha \\
                        \beta
                      \end{array}\right)_d=\left(\begin{array}{cccc}
\alpha_1,&\alpha_2,&\ldots,&\alpha_s\\[3pt]
\beta_1,&\beta_2,&\ldots,&\beta_t
\end{array} \right)_{d}
\end{equation}
be a partition in $V_2(1,n)$. By definition, we see that $\alpha_1=\beta_1=d$, $s(\beta)=\beta_t=1$ and $\ell(\alpha)-\ell(\beta)=s-t\le -2$.

Define
\begin{equation}
\psi_2(\lambda)=\left(\begin{array}{c}
                        \gamma \\
                        \delta
                      \end{array}\right)_{d'}=\left(\begin{array}{cccc}
\alpha_2,&\ldots,&\alpha_s\\[3pt]
\beta_2,&\ldots,&\beta_{t-1}
\end{array} \right)_{d+1}.
\end{equation}
It is evident that $\ell(\gamma)-\ell(\delta)=(s-1)-(t-2)=s-t+1\le -1$, $\gamma_1=\alpha_2\le\alpha_1= d=d'-1$ and $\delta_1=\beta_2\le \beta_1=d=d'-1$. Moreover,
\begin{align*}
|\psi_2(\lambda)| &=|\gamma|+|\delta|+(d')^2 \\
   &=\left(\sum_{i=2}^{s}\alpha_i\right)
   +\left(\sum_{i=2}^{t-1}\beta_i\right)+(d+1)^2 \\
  & = \left(|\alpha|-d\right)
   +\left(|\beta|-d-1\right)+(d+1)^2\\
  &=|\alpha|+|\beta|+d^2=|\lambda|.
\end{align*}
Thus $\psi_2(\lambda)\in U_2(0,n)$.

To show $\psi_2$ is a bijection,  we aim to build the inverse map $\psi_2^{-1}$ from $U_2(0,n)$ to $V_2(1,n)$.
Let
\begin{equation*}
  \mu=\left(\begin{array}{c}
              \gamma \\
              \delta
            \end{array}\right)_{d'}=\left(\begin{array}{cccc}
\gamma_1,&\ldots,&\gamma_{s'}\\[3pt]
\delta_1,&\ldots,&\delta_{t'}
\end{array} \right)_{d'}
  \end{equation*}
  be a partition in $U_2(0,n)$. By definition, we see that $\gamma_1\le d'-1$, $\delta_1\le d'-1$ and $s'-t'\le -1$. Define $\psi_2^{-1}(\mu)$ as follows:
  \begin{equation*}
  \psi_2^{-1}(\mu)=\left(\begin{array}{c}
                        \alpha \\
                        \beta
                      \end{array}\right)_d
  =\left(\begin{array}{ccccccc}
                        d'-1,&\gamma_1,&\ldots,&\gamma_{s'}& \\
                        d'-1,&\delta_1,&\ldots,&\delta_{t'},&1
                      \end{array}\right)_{d'-1}.
  \end{equation*}
It is easy to check that $\psi_2^{-1}(\mu)\in V_2(1,n)$ and $\psi_2^{-1}$ is the inverse map of $\psi_2$. This completes the proof.\qed

For example, let
$$\lambda=\left(\begin{array}{ccccccccccc}
                  5, & 4, &4,& 3, & 2, & 2 \\
                  5,  & 3, & 2, & 2,&2,&1,&1,&1,&1
                \end{array}\right)_{5}\in V_2(1,63).$$
Applying $\psi_2$ on $\lambda$, we obtain that
$$\psi_2(\lambda)=\mu=\left(\begin{array}{ccccccccccc}
                   4, &4,& 3, & 2, & 2 \\
                   3, & 2, & 2,&2,&1,&1,&1
                \end{array}\right)_{6}$$
which is in $U_2(0,63)$. Applying $\psi_2^{-1}$ on $\mu$, we recover $\lambda$.

We next establish the injection $\psi_3$ from $V_3(1,n)$ into $U_3(0,n)$.

\begin{lem}\label{lem-psi3}
  There is an injection $\psi_3$ from $V_3(1,n)$ into $U_3(0,n)$.
\end{lem}

\pf Let
\begin{equation}
\lambda=\left(\begin{array}{c}
                        \alpha \\
                        \beta
                      \end{array}\right)_d=\left(\begin{array}{cccc}
\alpha_1,&\alpha_2,&\ldots,&\alpha_s\\[3pt]
\beta_1,&\beta_2,&\ldots,&\beta_t
\end{array} \right)_{d}
\end{equation}
be a partition in $V_3(1,n)$.  By definition, we see  $\ell(\alpha)-\ell(\beta)=s-t\le -2$, $\alpha_1=\beta_1=d$ and $\beta_t=s(\beta)\ge 2$. Define $i=\max\{p\colon 1\le p\le s, \alpha_p\ge \beta_{p+2}-1\}$. From $\alpha_1=d> d-1=\beta_1-1\ge \beta_3-1$, we see that  such $i$ exists. The map $\psi_3(\lambda)$ is defined as follows:
\begin{equation*}
\psi_3(\lambda)=\left(\begin{array}{c}
                        \gamma \\
                        \delta
                      \end{array}\right)_{d'}
=\left(\begin{array}{ccccccccccc}
\alpha_2,&\ldots,&\alpha_i,&\beta_{i+2}-1,&\ldots,&\beta_t-1&\\[3pt]
\beta_2,&\ldots,&\beta_{i+1},&\alpha_{i+1}+1,&\ldots,&\alpha_s+1,&1^{t-s-2}
\end{array} \right)_{d+1},
\end{equation*}
where $1^{t-s-2}$ means $1$ appears $t-s-2$ times. We first show that $\gamma$ and $\delta$ are partitions. By the definition of $i$, we have $\alpha_i\ge \beta_{i+2}-1$, thus $\gamma$ is a partition. From the maximum of $i$, we find that either $i=s$ or both $i<s$ and $\alpha_{i+1}<\beta_{i+3}-1$.
If $i=s$, then $\delta=(\beta_2,\ldots,\beta_{s+1},1^{t-s-2})$. By $s-t\le -2$ we see that $\beta_{s+1}\ge\beta_t\ge 2$, thus $\delta$ is a partition. If $i<s$ and $\alpha_{i+1}<\beta_{i+3}-1$, then $\alpha_{i+1}+1<\beta_{i+3}\le \beta_{i+1}$. So in either case, $\delta$ is a partition. Thus $\psi_3$ is well defined.

We next verify $\psi_3(\lambda)\in U_3(0,n)$.
It is clear that $\gamma_1=\alpha_2\le \alpha_1=d=d'-1$, $\delta_1=\beta_2\le\beta_1= d=d'-1$ and
$\ell(\gamma)-\ell(\delta)=(t-2)-(t-2)=0.$
Moreover,
\begin{align*}
  |\psi_3(\lambda)| &=|\gamma|+|\delta|+(d')^2 \\
   &=\left(\sum_{k=2}^{i}\alpha_k\right)
   +\left(\sum_{k=i+2}^{t}\beta_k-1\right)+
  \left(\sum_{k=2}^{i+1}\beta_k\right)
   +\left(\sum_{k=i+1}^{s}\alpha_k+1\right)+t-s-2+(d+1)^2 \\
   &=|\alpha|-\alpha_1+|\beta|-\beta_1-(t-i-1)+(s-i)+t-s-2+(d+1)^2 \\
   &=|\alpha|+|\beta|+d^2=|\lambda|.
\end{align*}
Hence $\psi_3(\lambda)\in U_3(0,n)$.

We proceed to show that $\psi_3$ is an injection. To this end, let $H(0,n)\subseteq U_3(0,n)$ be the image set of $\psi_3$. Given
$$\mu=\left(\begin{array}{c}
           \gamma \\
           \delta
         \end{array}\right)_{d'}
=\left(\begin{array}{cccc}
\gamma_1, & \ldots,&\gamma_{t'} \\
\delta_1,& \ldots, &\delta_{t'} \end{array}\right)_{d'}\in H(0,n).$$
Let $j=\max\{p\colon 1\le p\le t', \delta_p\ge \gamma_p+1\}$.  By the construction of $\psi_3$, we see that $\delta_i=\beta_{i+1}$ and $\gamma_i=\beta_{i+2}-1$. Thus $\delta_i\ge \gamma_i+1$, which implies the existence of $j$.  By the definition of $U_3(0,n)$, we see that $\gamma_1\le d'-1$ and $\delta_1\le d'-1$.
We define a map $\pi$ from $H(0,n)$ to $V_3(1,n)$ as given below.
\[\pi(\mu)=\left(\begin{array}{c}
\alpha \\
\beta
\end{array}\right)_d=
\left(\begin{array}{ccccccccccccc}
        d'-1,& \gamma_1, & \ldots, & \gamma_{j-1},&\delta_{j+1}-1, &\ldots,&\delta_{t'}-1 \\
        d'-1, & \delta_1, & \ldots, & \delta_j,&\gamma_j+1,&\ldots,&\gamma_{t'}+1
\end{array}\right)_{d'-1}.\]
We first show that $\alpha$ and $\beta$ are partitions. By the definition of $j$, we have $\delta_j\ge \gamma_j+1$, which implies $\beta$ is a partition. From the maximum of $j$, we see that either $j=t'$ or both $j<t'$ and $\delta_{j+1}<\gamma_{j+1}+1$. If $j=t'$, then $\alpha=(d'-1, \gamma_1,  \ldots,  \gamma_{t'-1})$ which is clearly a partition. If $j<t'$, then $\delta_{j+1}-1<\gamma_{j+1}\le\gamma_{j-1}$. So in either case, $\alpha$ is a partition.

To show $\pi(\mu)\in V_3(1,n)$, it is enough to verify $\alpha_1=\beta_1=d$, $\ell(\alpha)-\ell(\beta)\le -2$, $s(\beta)\ge 2$ and $|\pi(\mu)|=|\mu|$. On the one hand, it is clear that $\alpha_1=\beta_1=d'-1=d$. Moreover $\ell(\alpha)\le t'$, $\ell(\beta)=t'+2$, which implies $\ell(\alpha)-\ell(\beta)\le -2$. Furthermore, $s(\beta)=\gamma_{t'}+1\ge 2$. On the other hand,
\begin{align*}
  |\pi(\mu)| &=|\alpha|+|\beta|+d^2 \\
   &=2d'-2+\left(\sum_{i=1}^{j-1}\gamma_i\right)
   +\left(\sum_{i=j+1}^{t'}\delta_i-1\right)
   +\left(\sum_{i=1}^{j}\delta_i\right)
   +\left(\sum_{i=j}^{t'}\gamma_i+1\right)+(d'-1)^2\\
   &=|\gamma|+|\delta|-(t'-j)+(t'-j+1)+{d'}^2-1 \\
  &=|\gamma|+|\delta|+{d'}^2=|\mu|.
\end{align*}
Hence $\pi(\mu)\in V_3(1,n)$.

We conclude this proof by showing that for any $\lambda\in V_3(1,n)$, $
  \pi(\psi_3(\lambda))=\lambda.$
From the construction of $\psi_3$ and $\pi$, it sufficient to show the integer $i$ appears in $\psi_3(\lambda)$ coincides with the integer $j$ appears in $\pi(\psi_3(\lambda))$.

Recall that $i$ is the maximum integer $1\le i\le s$ such that $\alpha_i\ge \beta_{i+2}-1$. Hence either $i=s$ or both $i<s$ and $\alpha_{k}<\beta_{k+2}-1$ for all $s\ge k\ge i+1$. If $i=s$, then $$\gamma=(\alpha_2,\ldots,\alpha_s,\beta_{s+2}-1,\ldots,\beta_t-1)$$
and
\[\delta=(\beta_2,\ldots,\beta_{s+1},1^{t-s-2}).\]
Clearly $\gamma_s=\beta_{s+2}-1$ and $\delta_s=\beta_{s+1}$, which implies $\gamma_s\le \delta_s-1$, thus $j\ge s$. Suppose that $j> s$, then $\gamma_j=\beta_{j+2}-1\ge \beta_t-1\ge 1$ and $\delta_j=1$. Hence $\gamma_j>\delta_j-1$, contradict to the definition of $j$.   So in this case we have $j=s=i$.

If $i<s$ and for any $t-2\ge k\ge i+1$, we have $\alpha_{k}<\beta_{k+2}-1$. Then $\gamma_i=\beta_{i+2}-1$ and $\delta_i=\beta_{i+1}$, which implies $\gamma_i\le \delta_i-1$. Hence $j\ge i$. Moreover, for any $t-2\ge k\ge i+1$, by definition $\gamma_k=\beta_{k+2}-1$ and $\delta_k=\alpha_k+1$. Thus $\gamma_k=\beta_{k+2}-1>\alpha_k=\delta_k-1$. This yields $j\le i$.
So in either case, $i=j$. This completes the proof.\qed

For instance, let
$$\lambda=\left(\begin{array}{ccccccccccc}
                  7, & 6, &5,& 3, & 1 \\
                  7,  &7,& 6, & 6, & 4,&3,&3,&2
                \end{array}\right)_{7}\in V_3(1,109).$$
Applying $\psi_3$ on $\lambda$, we see that $i=4$, hence
$$\psi_3(\lambda)=\mu=\left(\begin{array}{ccccccccccc}
                   6, &5,& 3, & 2, & 2,&1 \\
                   7, & 6, & 6,&4,&2,&1
                \end{array}\right)_{8}$$
which is in $U_3(0,109)$. Applying $\pi$ on $\mu$,  we have $j=4$ and $\pi(\mu)=\lambda$.

We are now in a position to prove Theorem \ref{thm-main} for $m=0$.

{\noindent \it Proof of Theorem \ref{thm-main} for $m=0$.} Let $\lambda$ be a partition in $V(1,n)$. Using Lemma \ref{lem-psi2} and Lemma \ref{lem-psi3}, we may define the injection $\Psi$ from $V(1,n)$ into $U(0,n)$ as follows:
\begin{equation}\label{equ-def-Psi}
  \Psi(\lambda)=\begin{cases}
                  \lambda, & \mbox{if } \lambda\in V_1(1,n); \\
                  \psi_2(\lambda), & \mbox{if } \lambda\in V_2(1,n); \\
                  \psi_3(\lambda), & \mbox{if } \lambda\in V_3(1,n).
                \end{cases}
\end{equation}

We next show that $u(0,n)>u(1,n)$ for $n\ge 6$. From the construction of $\Psi$, we find that the set $U_4(0,n)$ and $U_5(0,n)$  contain no image of $\Psi$, which implies
\begin{equation}\label{ine-est-u0n-u1n-u4n-u5n}
  u(0,n)-u(1,n)\ge \# U_4(0,n)+\# U_5(0,n).
\end{equation}
We proceed to show that $U_4(0,n)\cup U_5(0,n)\ne \emptyset$. There are two cases.

Case 1: $n=2t$ for $t\ge 3$. Define
\[\mu=\left(\begin{array}{cc}
                  1^{t-3} & \\
                  2, & 1^{t-3}
                \end{array}\right)_2,\]
which is clearly in $U_4(0,n)$.

Case 2: $n =2t+1$ for $t\ge 3$. Set
\[\mu=\left(\begin{array}{cc}
                  1^{t-2} & \\
                  2, & 1^{t-3}
                \end{array}\right)_2,\]
which is clearly in $U_5(0,n)$. This completes the proof.\qed

\section{Connection to $\ospt(n)$}\label{5}

This section is devoted to prove Theorem \ref{cor-ospt-lower} and Theorem \ref{thm-asy-ospt-pn4}. To this end, we first give a new combinatorial interpretation of $u(1,n)$, namely Lemma \ref{lem-com-int-u1n}. Then we establish a connection between $u(0,n)-u(1,n)$ and $\ospt(n)$, that is Corollary \ref{cor-u0n-ospt}. From Corollary \ref{cor-u0n-ospt}, we find that  Theorem \ref{thm-asy-ospt-pn4} can be deduced from \eqref{eq-asy-umn} and \eqref{Mao-asy}. Moreover, by Corollary \ref{cor-u0n-ospt}, Theorem \ref{cor-ospt-lower} is equivalent to a lower bound of $u(0,n)-u(1,n)$, which is  Theorem \ref{thm-est-u0n-u1n}. Finally, we prove  Theorem \ref{thm-est-u0n-u1n}. This completes the proof of Theorem \ref{cor-ospt-lower}.

We first present a new combinatorial interpretation of $u(1,n)$.

\begin{lem}\label{lem-com-int-u1n}
  For $n\ge 1$, $u(1,n)$ counts the number of partitions $\lambda$ of $n$ such that ${\rm rank}(\lambda)\le 0$ and the rank-set of $\lambda$ is not contain $-1$.
\end{lem}

\pf Let $A(1,n)$ denote the set of partitions $\lambda$ of $n$ such that ${\rm rank}(\lambda)\le 0$ and   $-1$ is not in the rank-set of $\lambda$. We use $B(1,n)$ to denote the set of partitions $\mu$ of $n$ such that ${\rm rank}(\lambda)\le -2$ and $0$ is in the rank-set of $\lambda$. By Theorem \ref{thm-com-int-umn}, we have $\#B(1,n)=u(1,n)$. Thus in order to prove Lemma \ref{lem-com-int-u1n}, it sufficient to build a bijection $\rho$ between $A(1,n)$ and $B(1,n)$.

Given $\lambda\in A(1,n)$.  Since $-1$ is not contained in the rank-set of $\lambda$, let $j$ be the maximum integer such that $j-\lambda_{j+1}<-1$. Thus $\lambda_{j+1}\ge  j+2$. Recall that the rank-set is a strictly increasing integer sequence. From the maximum of $j$, we see that $j+1-\lambda_{j+2}>-1$, which implies $j+1\ge \lambda_{j+2}$. Define
\begin{equation}\label{equ-def-rho}
  \mu=\rho(\lambda)=(\lambda_1-1,\ldots,\lambda_{j+1}-1,j+1,\lambda_{j+2},\ldots,\lambda_{\ell(\lambda)}).
\end{equation}
From the above analysis, we see that $\rho(\lambda)$ is a partition and $\mu_{j+2}=j+1$, which implies $j+1-\mu_{j+2}=0$ is in the rank-set of $\mu$. Moreover,
\begin{equation*}
  {\rm rank}(\mu)=
  (\lambda_1-1)-(\ell(\lambda)+1)={\rm rank}(\lambda)-2\le -2.
\end{equation*}
It is trivial to check that $|\mu|=|\lambda|$. Hence $\mu\in B(1,n)$.

We proceed to construct the inverse map $\rho^{-1}$ of $\rho$. For any $\mu\in B(1,n)$, since $0$ is in the rank-set of $\mu$, there exists an integer $k\ge 1$ such that $k-\mu_{k+1}=0$.  Define
\begin{equation}\label{equ-def-rho-1}
 \rho^{-1}(\mu)=\lambda=(\mu_1+1,\ldots,\mu_k+1,\mu_{k+2},\ldots,\mu_{\ell(\mu)}).
\end{equation}
By $\mu_{k+1}=k$, we see that $\mu_k\ge k\ge \mu_{k+2}$, thus
\begin{equation}\label{eq-te-1}
  (k-1)-\lambda_k=k-1-\mu_k-1\le -2,
\end{equation}
and
\begin{equation}\label{eq-te-2}
k-\lambda_{k+1}=k-\mu_{k+2}\ge 0.
\end{equation}
Since the rank-set of $\lambda$ is a strict increasing sequence, by \eqref{eq-te-1} and \eqref{eq-te-2} we see that $-1$ is not contained in the rank-set of $\lambda$. Moreover ${\rm rank}(\lambda)=(\mu_1+1)-(\ell(\mu)-1)={\rm rank}(\mu)+2\le 0.$ It is trivial to check that  $|\lambda|=|\mu|$.
Thus $\rho^{-1}(\mu)\in A(1,n)$.

To show $\rho^{-1}$ is the inverse map of $\rho$, from the construction of $\rho$ and $\rho^{-1}$, it sufficient to show that the integer $k$ in \eqref{equ-def-rho-1} equals the integer $j$ in \eqref{equ-def-rho} plus $1$, that is $k=j+1$.
Actually, this equation follows from \eqref{eq-te-1}, \eqref{eq-te-2} and the definition of $j$. This completes the proof.\qed

From Lemma \ref{lem-com-int-u1n},
we find the following relation between $\ospt(n)$ and $u(0,n)-u(1,n)$, which plays a crucial role in the proof of Theorem \ref{cor-ospt-lower} and Theorem \ref{thm-asy-ospt-pn4}.

\begin{core}\label{cor-u0n-ospt}
For $n\ge 0$,
\begin{equation}\label{equ-co-u0n}
  \ospt(n)=\frac{p(n)}{4}+\frac{N(0,n)}{4}+\frac{u(0,n)-u(1,n)}{2}.
\end{equation}
\end{core}

\pf  From Theorem \ref{thm-com-int-umn}, we have $u(0,n)$ equals the number   of partitions $\lambda$ of $n$ such that ${\rm rank}(\lambda)\le 0$ and $-1$ is in the rank-set of $\lambda$. By Lemma \ref{lem-com-int-u1n},  $u(1,n)$ counts the number of partitions $\lambda$ of $n$ such that ${\rm rank}(\lambda)\le 0$ and $-1$ is not in the rank-set of $\lambda$. Thus $u(0,n)+u(1,n)$ is equal to the   number of partitions $\lambda$ of $n$ such that ${\rm rank}(\lambda)\le 0$. In other words,
\begin{equation}\label{equ-u0n+u1n}
  u(0,n)+u(1,n)=\sum_{m=-\infty}^{0}N(m,n).
\end{equation}
Using
$
 N(m,n)=N(-m,n)
$
due to Dyson \cite{Dyson-1944}, we deduce that
\begin{align}\label{equ-pn-n0n}
  p(n)+N(0,n)&= \sum_{m=-\infty}^{+\infty}N(m,n)+N(0,n)\nonumber\\
  &=\left(\sum_{m=-\infty}^{-1}N(m,n)\right)
  +\left(\sum_{m=1}^{\infty}N(m,n)\right)+2N(0,n)\nonumber\\
  &=2\left(\sum_{m=-\infty}^{-1}N(m,n)+N(0,n)\right)\nonumber\\
  &=2\sum_{m=-\infty}^{0}N(m,n).
\end{align}
From \eqref{equ-u0n+u1n} and \eqref{equ-pn-n0n}, we see that
\begin{equation*}
  u(0,n)= \frac{1}{2}\left(\sum_{m=-\infty}^{0}N(m,n)+u(0,n)-u(1,n)\right)
  =\frac{p(n)}{4}+\frac{N(0,n)}{4}+\frac{u(0,n)-u(1,n)}{2}.
\end{equation*}
By \eqref{equ-cor-ospt-u0n}, we derive \eqref{equ-co-u0n}.
\qed

We next show that Theorem \ref{thm-asy-ospt-pn4} is a direct consequence of Corollary \ref{cor-u0n-ospt}.

{\noindent\it Proof of Theorem \ref{thm-asy-ospt-pn4}.} Set  $m=0$ in \eqref{eq-asy-umn}, we deduce that
\[u(0,n)-u(1,n)\sim \frac{\pi}{96\sqrt{2}n^{\frac{3}{2}}}e^{\pi\sqrt{\frac{2n}{3}}}.\]
Combining with \eqref{Mao-asy} and Corollary \ref{cor-u0n-ospt}, we deduce \eqref{ospt-ays}.\qed

Using Corollary \ref{cor-u0n-ospt}, it is clear that Theorem \ref{cor-ospt-lower} is equivalent to the following theorem.

\begin{thm}\label{thm-est-u0n-u1n}
  For $n\ge 6$, we have
  \begin{equation}\label{ine-est-u0n-u1n}
    u(0,n)-u(1,n)> \frac{N(0,n)}{4}.
  \end{equation}
\end{thm}

  The rest part of this section is devoted to give a proof of Theorem \ref{thm-est-u0n-u1n}.
By \eqref{ine-est-u0n-u1n-u4n-u5n}, Theorem \ref{thm-est-u0n-u1n} can be deduced from the following inequality:
\begin{equation}\label{ine-u4n-u5n-n0n}
  4\# U_4(0,n)+4\# U_5(0,n)> N(0,n).
\end{equation}
We now prove \eqref{ine-u4n-u5n-n0n}. Let $P(0,n)$ denote the set of partitions of $n$ with rank $0$. Clearly $\#P(0,n)=N(0,n)$. Given  $(\alpha,\beta)_d\in P(0,n)$, by definition, we have $\ell(\alpha)=\ell(\beta)$, $\alpha_1\le d$ and $\beta_1\le d$. In order to show \eqref{ine-u4n-u5n-n0n},  we further divide $P(0,n)$ into ten disjoint subsets $P_i(0,n)$ $(1\le i\le 10)$ as follows:
\begin{itemize}
  \item[(1)] $P_1(0,n)$ is the subset of Durfee symbols $(\alpha,\beta)_d$ in $P(0,n)$ such that $\beta_1=d$ and $\alpha_1<d$;
  \item[(2)] $P_2(0,n)$ is the subset of Durfee symbols $(\alpha,\beta)_d$ in $P(0,n)$ such that $\beta_1<d$ and $\alpha_1=d$;
  \item[(3)] $P_3(0,n)$ is the subset of Durfee symbols $(\alpha,\beta)_d$ in $P(0,n)$ such that $\beta_1=d$ and $\alpha_1=d>\alpha_2$;
  \item[(4)] $P_4(0,n)$ is the subset of Durfee symbols $(\alpha,\beta)_d$ in $P(0,n)$ such that $\beta_1=d$ and $\alpha_1=\alpha_2=d>\alpha_3$;
  \item[(5)] $P_5(0,n)$ is the subset of Durfee symbols $(\alpha,\beta)_d$ in $P(0,n)$ such that $d\ne 2$, $\beta_1=d$ and $\alpha_1=\alpha_2=\alpha_3=d$;
  \item[(6)]$P_6(0,n)$ is the subset of Durfee symbols $(\alpha,\beta)_d$ in $P(0,n)$ such that $d=2$, $\beta_1=\beta_2=2$ and $\alpha_1=\alpha_2=\alpha_3=2$;
  \item[(7)] $P_7(0,n)$ is the subset of Durfee symbols $(\alpha,\beta)_d$ in $P(0,n)$ such that $d=2$, $\beta_1=2>\beta_2$ and $\alpha_1=\alpha_2=\alpha_3=2$;
  \item[(8)] $P_8(0,n)$ is the subset of Durfee symbols $(\alpha,\beta)_d$ in $P(0,n)$  such that $\beta_1<d$ and $\alpha_1\le d-2$;
  \item[(9)] $P_9(0,n)$ is the subset of Durfee symbols $(\alpha,\beta)_d$ in $P(0,n)$ such that  $\beta_1<d$ and $\alpha_1=d-1>\alpha_2$;
  \item[(10)] $P_{10}(0,n)$ is the subset of Durfee symbols $(\alpha,\beta)_d$ in $P(0,n)$ such that  $\beta_1<d$ and $\alpha_1=d-1=\alpha_2$.
\end{itemize}

We then build eight injections $\chi_i$ $(1\le i\le 8)$. To be specific,
\begin{align*}
  \chi_1 \colon & P_1(0,n)\rightarrow U_5(0,n);  \\
  \chi_2 \colon & P_2(0,n)\rightarrow U_5(0,n); \\
  \chi_3 \colon & P_3(0,n)\rightarrow U_4(0,n);\\
  \chi_4 \colon & P_4(0,n)\cup P_7(0,n)\rightarrow U_4(0,n); \\
  \chi_5 \colon & P_5(0,n)\cup P_6(0,n)\rightarrow U_5(0,n); \\
  \chi_6 \colon & P_8(0,n)\rightarrow U_4(0,n); \\
  \chi_7 \colon & P_9(0,n)\rightarrow U_4(0,n);\\
  \chi_8 \colon & P_{10}(0,n)\rightarrow U_5(0,n).
\end{align*}

It is trivial to check that $U_5(0,n)$ is coincide with $P_1(0,n)$. So we may set the map $\chi_1$ to be the identity map. The following lemma describe the bijection $\chi_2$.

\begin{lem}\label{lem-chi2}
  There is a bijection $\chi_2$ between $P_2(0,n)$ and $U_5(0,n)$.
\end{lem}

\pf Let $\lambda=(\alpha,\beta)_d\in P_2(0,n)$. By definition, $\ell(\alpha)=\ell(\beta)$, $\beta_1<d$ and $\alpha_1=d$. Define
\[\chi_2(\lambda)=\left(\begin{array}{c}
                    \beta \\
                    \alpha
                  \end{array}\right)_d.\]
It is trivial to check that $\chi_2(\lambda)\in U_5(0,n)$ and $\chi_2$ is a bijection. \qed

We next present the injection $\chi_3$.

\begin{lem}
  For $n\ge 4$, there is an injection $\chi_3$ from $P_3(0,n)$ into $U_4(0,n)$.
\end{lem}

\pf We first show that if $(\alpha,\beta)_d\in P_3(0,n)$, then $d\ge 2$. Assume the contrary, if $d=1$, then $\alpha_1=1>\alpha_2$ which yields $\alpha=(1)$. Since $\ell(\alpha)=\ell(\beta)$, we see that $\beta=(1)$. Thus $|\lambda|=|\alpha|+|\beta|+d^2=3$, contradicts to $n\ge 4$. Hence $d\ge 2$.

By definition, we have $\beta_1=d$, $\alpha_1=d>\alpha_2$ and $\ell(\alpha)=\ell(\beta)$. Define
\begin{equation}\label{equ-def-chi3}
  \chi_3(\lambda)=\left(\begin{array}{c}
                    \gamma \\
                    \delta
                  \end{array}\right)_{d'}
  =\left(\begin{array}{ccccc}
  \alpha_1-1,& \alpha_2,& \ldots, & \alpha_{\ell(\alpha)} \\
  \beta_1,&\ldots, & \beta_{\ell(\beta)}, & 1
  \end{array}\right)_d.
\end{equation}
It is clear that $\ell(\gamma)-\ell(\delta)=\ell(\alpha)-\ell(\beta)-1=-1$, $\gamma_1=\alpha_1-1=d-1<d=d'$ and $\delta_1=\beta_1=d=d'$. Moreover, $|\lambda|=|\chi_3(\lambda)|$, thus $\chi_3(\lambda)\in U_4(0,n)$.

To show $\chi_3$ is an injection, let $K_3(0,n)$ be the image set of $\chi_3$, which has been shown to be a subset of $U_4(0,n)$. For any $\mu=(\gamma,\delta)_{d'}\in K_3(0,n)$, from the construction of $\chi_3$, we see that $\gamma_1=\alpha_1-1=d-1=d'-1$ and $s(\delta)=1$. Define the map $\zeta_3$ as follows:
\[\zeta_3(\mu)=\left(\begin{array}{c}
                    \alpha \\
                    \beta
                  \end{array}\right)_{d}
  =\left(\begin{array}{ccccc}
  \gamma_1+1,& \gamma_2,& \ldots, & \gamma_{\ell(\gamma)} \\
  \delta_1,&\ldots, & \delta_{\ell(\delta)-1}
  \end{array}\right)_{d'}\]
It is trivial to check that for any $\lambda\in P_3(0,n)$,  $\zeta_3(\chi_3(\lambda))=\lambda$. Thus $\chi_3$ is an injection. \qed

For instance, let
$$\lambda=\left(\begin{array}{ccccccccccc}
                  5, & 4, &4,& 2, & 2 \\
                  5,  &5,& 3, & 3, & 1
                \end{array}\right)_{5}\in P_3(0,59).$$
Applying $\chi_3$ on $\lambda$,
$$\chi_3(\lambda)=\mu=\left(\begin{array}{ccccccccccc}
                  4, & 4, &4,& 2, & 2 \\
                  5,  &5,& 3, & 3, & 1,&1
                \end{array}\right)_{5}$$
which is in $U_4(0,59)$. Applying $\zeta_3$ on $\mu$, we recover $\lambda$.

We next establish the injection $\chi_4$.

\begin{lem}
  For $n\ge 6$, there is an injection $\chi_4$ from $P_4(0,n)\cup P_7(0,n)$ into $U_4(0,n)$.
\end{lem}

\pf We divide the set $U_4(0,n)$ into two disjoint subsets as follows:
\begin{itemize}
  \item[(1)] $U_4^1(0,n)$ denote the subset $(\gamma,\delta)_{d'}$ of $U_4(0,n)$ such that $2$ is a part of $\delta$;
  \item[(2)] $U_4^2(0,n)$ denote the subset $(\gamma,\delta)_{d'}$ of $U_4(0,n)$ such that $2$ is not a part of $\delta$.
\end{itemize}
We next build an injection $\eta_1\colon P_4(0,n)\rightarrow U_4^1(0,n)$ and an injection $\eta_2\colon P_7(0,n)\rightarrow U_4^2(0,n)$. For any $\lambda\in P_4(0,n)\cup P_7(0,n)$, the injection $\chi_4$ can be defined as follows:
\begin{equation}\label{eq-def-chi4}
  \chi_4(\lambda)=\begin{cases}
                    \eta_1(\lambda), & \mbox{if } \lambda\in P_4(0,n); \\
                    \eta_2(\lambda), & \mbox{if } \lambda\in P_7(0,n).
                  \end{cases}
\end{equation}

We first describe the  injection $\eta_1$. Given $\lambda=(\alpha,\beta)_d\in P_4(0,n)$, by definition, we have $\ell(\alpha)=\ell(\beta)$ and $\alpha_1=\alpha_2=\beta_1=d>\alpha_3$. We first show that $d\ge 2$. Assume the contrary, if $d=1$, then by definition we have $\alpha=(1,1)$. From $\ell(\alpha)=\ell(\beta)$ we see that $\beta=(1,1)$. Thus $|\lambda|=|\alpha|+|\beta|+d^2=5$, contradicts to $n\ge 6$. Thus $d\ge 2$.
For $(\alpha,\beta)_d\in P_4(0,n)$, let $k$ be the minimum integer such that $\beta_k=1$ (if such $k$ not exists, set $k=\ell(\beta)+1$). Since $\beta_1=d\ge 2$, we see that $k\ge 2$. Define
\begin{equation}\label{equ-def-chi4-1}
  \eta_1(\lambda)=\left(\begin{array}{c}
                    \gamma \\
                    \delta
                  \end{array}\right)_{d'}
  =\left(\begin{array}{ccccccc}
  \alpha_1-1,& \alpha_2-1,&\alpha_3,& \ldots, & \alpha_{\ell(\alpha)}& \\
  \beta_1,&\ldots,& \beta_{k-1},&2,&\beta_k,& \ldots,& \beta_{\ell(\beta)}
  \end{array}\right)_d.
\end{equation}
From $\alpha_1=\alpha_2>\alpha_3$, we see that $\gamma$ is a partition. By  the choice of $k$, we see that $\beta_{k-1}\ge 2>\beta_k$. Hence $\delta$ is a partition. Thus the map $\eta_1$ is well defined.
 From the construction of $\eta_1$, we see that $\ell(\gamma)-\ell(\delta)=\ell(\alpha)-(\ell(\beta)+1)=-1,$
$\gamma_1=\alpha_1-1=d-1=d'-1$ and $\delta_1=\beta_1=d=d'$. Moreover, $2$ is a part of $\delta$ and $|\lambda|=|\eta_1(\lambda)|$.
Therefore, $\eta_1(\lambda)\in U_4^1(0,n)$.

We next show that $\eta_1$ is an injection. Let $K_{4}^1(0,n)$ denote the image set of $\eta_1$, which has been shown to be a subset of $U_4^1(0,n)$. For any $\mu=(\gamma, \delta)_{d'}\in K_{4}^1(0,n)$, from the construction of $\eta_1$, we see that $\gamma_1=\gamma_2=d'-1$ and there exists an integer $k$ such that $\delta_k=2$. Define
\[\zeta_{4}^1(\mu)=\left(\begin{array}{cccccccccc}
\gamma_1+1, & \gamma_2+1, & \gamma_3,& \ldots, & \gamma_{\ell(\gamma)} \\
\delta_1, & \ldots, & \delta_{k-1}, & \delta_{k+1}, & \ldots,&\delta_{\ell(\delta)}
                      \end{array}\right)_{d'}.\]
It is routine to check that $\zeta_{4}^1(\mu)\in P_4(0,n)$ and for any $\lambda\in P_4(0,n)$, $\zeta_{4}^1(\eta_1(\lambda))=\lambda$. Thus $\eta_1$ is an injection.

We next build an injection $\eta_2$ from $P_7(0,n)$ into the set $U_4^2(0,n)$. Given $\lambda=(\alpha,\beta)_d\in P_7(0,n)$, by definition  $d=2$ and $\alpha_1=\alpha_2=\alpha_3=\beta_1=2>\beta_2$. Assume $\ell(\alpha)=\ell(\beta)=t$, since $\alpha_3=2>\beta_2$, we see that $t\ge 3$ and $\beta=(2,1^{t-1})$. Define
\begin{equation}\label{eta2}
  \eta_2(\lambda)=\left(\begin{array}{c}
                          \gamma \\
                          \delta
                        \end{array}\right)_{d'}
  =\left(\begin{array}{ccccc}
           \alpha_4, & \ldots, & \alpha_t, & 1 \\
           3,&1^{t-2}  & &
         \end{array}\right)_3.
\end{equation}
It is clear that $2$ is not a part of $\delta$,  $\ell(\gamma)=t-2$ and $\ell(\delta)=t-1$, which implies $\ell(\gamma)-\ell(\delta)=-1$. Moreover $\delta_1=3=d'$ and $\gamma_1=\alpha_4\le 2<d'$. Furthermore,
\begin{align*}
  |\eta_2(\lambda)|&=\left(\sum_{i=4}^{t}\alpha_i\right)+1+3+t-2+(d')^2 \\
   &=|\alpha|-6+t+2+9 \\
  &=|\alpha|+(2+(t-1))+4 \\
  &=|\alpha|+|\beta|+4=|\lambda|.
\end{align*}
Hence $\eta_2(\lambda)\in U_4^2(0,n)$.

To show that $\eta_2$ is an injection, let $K_4^2(0,n)$ be the image set of $\eta_2$, which has been shown to be a subset of $U_4^2(0,n)$. Given $\mu=(\gamma,\delta)_{d'}\in K_4^2(0,n)$, let $\ell(\delta)=t'$, which implies $\ell(\gamma)=t'-1$. By the definition of $\eta_2$, we see that
$d'=3$, $\delta=(3,1^{t'-1})$ and $\gamma_{t'-1}=1$. Define the map $\zeta_4^2$ as follows:
\[\zeta_4^2(\mu)=\left(\begin{array}{c}
                      \alpha \\
                      \beta
                    \end{array}\right)_d=
\left(\begin{array}{cccccccccccc}
        2, & 2, & 2, & \gamma_1,&\ldots,&\gamma_{t'-2} \\
        2,& 1^{t'} & & &
      \end{array}\right)_2.\]
It is not difficult to check that $\zeta_4^2(\mu)\in P_7(0,n)$ and $\zeta_4^2(\eta_2(\lambda))=\lambda$ for any $\lambda\in P_7(0,n)$. Hence $\eta_2$ is an injection. This completes the proof. \qed

For example, let
$$\lambda=\left(\begin{array}{ccccccccccc}
                  6, & 6, &5,& 3, & 3,&2,&2,&1,&1 \\
                  6,  &6,& 6, & 5, & 4,&2,&1,&1,&1
                \end{array}\right)_{6}\in P_4(0,97).$$
Applying $\eta_1$ on $\lambda$, we see that $k=7$ and
$$\eta_1(\lambda)=\left(\begin{array}{c}
                          \gamma \\
                          \delta
                        \end{array}\right)_6=
\left(\begin{array}{ccccccccccc}
                  5, & 5, &5,& 3, & 3,&2,&2,&1,&1 \\
                  6,  &6,& 6, & 5, & 4,&2,&2,&1,&1,&1
\end{array}\right)_{6}$$
which is in $U_4^1(0,97)$. Applying $\zeta_{4}^1,$ we recover $\lambda.$

For another example, let
\[\lambda=\left(\begin{array}{ccccccc}
                  2, & 2, & 2, & 2, & 2, & 1, & 1\\
                  2, & 1, & 1, & 1, & 1, & 1, & 1
                \end{array}\right)_2,\]
which is in $P_7(0,24)$. We see that $t=7$ and
\[\eta_2(\lambda)=\mu=\left(\begin{array}{ccccccc}
                    2, & 2, & 1, & 1,&1\\
                  3, & 1, & 1, & 1, & 1, & 1
                \end{array}\right)_3\]
which is in $U_4^2(0,24)$. Applying $\zeta_4^2$ on $\mu$, we recover $\lambda.$

We now describe the injection $\chi_5$.

\begin{lem}
  For $n\ge 9$, there is an injection $\chi_5$ from $P_5(0,n)\cup P_6(0,n)$ into $U_5(0,n)$.
\end{lem}

\pf We divide $U_5(0,n)$ into two disjoint subsets $U_5^1(0,n)$ and $U_5^2(0,n)$ as given below.
\begin{itemize}
  \item[(1)] $U_5^1(0,n)$ denote the subset $(\gamma,\delta)_{d'}$ of $U_5(0,n)$ such that $d'\ne 3$;
  \item[(2)] $U_5^2(0,n)$ denote the subset $(\gamma,\delta)_{d'}$ of $U_5(0,n)$ such that $d'=3$.
\end{itemize}
We next build an injection $\kappa_1$ from $P_5(0,n)$ into $U_5^1(0,n)$ and an injection $\kappa_2$ from $P_6(0,n)$ into $U_5^2(0,n)$. Clearly,
the injection $\chi_5$ can be defined as follows:
\begin{equation*}
  \chi_5(\lambda)=\begin{cases}
                    \kappa_1(\lambda), & \mbox{if } \lambda\in P_5(0,n);\\
                    \kappa_2(\lambda), & \mbox{if }\lambda\in P_6(0,n).
                  \end{cases}
\end{equation*}

We first describe the injection $\kappa_1$. Given $\lambda=(\alpha,\beta)_d\in P_5(0,n)$, by definition $\alpha_1=\alpha_2=\alpha_3=\beta_1=d$ and either $d\ge 3$ or $d=1$. There are two cases.

Case 1: $d\ge 3$. In this case, let $\ell(\alpha)=\ell(\beta)=t\ge 3$. Denote $k$ to be the maximum integer $1\le k\le t$ such that $\alpha_k\ge d-2$. From $\alpha_3=d$ we see that $k\ge 3$. The map $\kappa_1$ can be defined as follows:
\begin{equation}\label{eta1-1}
  \kappa_1(\lambda)=\left(\begin{array}{c}
                          \gamma \\
                          \delta
                        \end{array}\right)_{d'}
  =\left(\begin{array}{cccccccccccc}
           \beta_2, & \ldots, & \beta_t &&&&&  \\
           d+1,&\alpha_4,  &\ldots, &\alpha_k,&d-2,&\alpha_{k+1},&\ldots,&\alpha_t.
         \end{array}\right)_{d+1}.
\end{equation}
From the choice of $k$, we see that $\alpha_k\ge d-2>\alpha_{k+1}$, hence $\delta$ is  a partition. Moreover   $d'=d+1\ge 4$, $\ell(\gamma)=\ell(\delta)=t-1$, $\delta_1=d+1=d'$ and $\gamma_1=\beta_2\le d=d'-1$. Furthermore,
\begin{align*}
  |\kappa_1(\lambda)| &=\left(\sum_{i=2}^{t}\beta_i\right)+(d+1)+(d-2)+
  \left(\sum_{i=4}^{t}\alpha_i\right)+(d+1)^2 \\
  &=|\beta|-d+|\alpha|-3d+(2d-1)+(d+1)^2 \\
   &=|\alpha|+|\beta|+d^2=|\lambda|.
\end{align*}
Hence, $\kappa_1(\lambda)\in U_5^1(0,n)$.

To show $\kappa_1$ in Case 1 is an injection, let $K_5^1(0,n)$ denote the image set of $\kappa_1$ in Case 1, which has been   shown to be a subset of $U_5^1(0,n)$. Given $\mu=(\gamma,\delta)_{d'}\in K_5^1(0,n)$, from the construction of $\kappa_1$, we see that $d'\ge 4$, $\delta_1=d'>\delta_2$ and there exists an integer $r$ such that $\delta_r=d-2=d'-3$. Define
\begin{align*}
 \zeta_5^1(\mu) & =\left(\begin{array}{c}
                      \alpha \\
                      \beta
                    \end{array}\right)_d \\
   & =
\left(\begin{array}{cccccccccccc}
        d'-1, & d'-1, & d'-1, & \delta_2,&\ldots,&\delta_{r-1},&\delta_{r+1},&
        \ldots,&\delta_{\ell(\delta)} \\
        d'-1,& \gamma_1, &\ldots, &\gamma_{\ell(\gamma)} &&&&&&
      \end{array}\right)_{d'-1}.
\end{align*}
It is trivial to check that $\zeta_5^1(\mu)\in P_5(0,n)$ and for any $\lambda=(\alpha,\beta)_d\in P_5(0,n)$ with $d\ge 3$, we have $\zeta_5^1(\kappa_1(\lambda))=\lambda$, which implies $\kappa_1$ is an injection in Case 1.

Case 2: $d=1$. In this case, since $\ell(\alpha)=\ell(\beta)$, we see that there is only one partition in $P_5(0,n)$, namely
\begin{equation*}
  \lambda=\left(\begin{array}{c}
                  1^t \\
                  1^t
                \end{array}\right)_1.
\end{equation*}
Thus $n=2t+1$, from $n\ge 9$ we have $t\ge 4$. Define
\begin{equation*}
  \kappa_1(\lambda)=\left(\begin{array}{ccc}
                    1^{t-2} &  \\
                    2, & 1^{t-3}
                  \end{array}\right)_2,
\end{equation*}
we see that
\[|\kappa_1(\lambda)|=2+(t-3)+(t-2)+4=2t+1=|\lambda|.\]
It is routine to check that $\kappa_1(\lambda)\in U_5^1(0,n)$.

Clearly $\kappa_1$ in Case 2 is also an injection, since there is only one partition in this case. Notice that the image set of $\kappa_1$ in the above two cases are not intersect. This fact follows from in Case 1 we have $d'\ge 4$ and in Case 2 we have $d'=2$. Hence we deduce that $\kappa_1$ is an injection.

We next establish the injection $\kappa_2$. Given $\lambda=(\alpha,\beta)_d\in P_6(0,n)$, by definition, we see that $\alpha_1=\alpha_2=\alpha_3=\beta_1=\beta_2=d=2$ and $\ell(\alpha)=\ell(\beta)$. Define
\begin{equation}\label{zeta2}
  \kappa_2(\lambda)=\left(\begin{array}{c}
                          \gamma \\
                          \delta
                        \end{array}\right)_{d'}
  =\left(\begin{array}{cccccccccccc}
           \alpha_4, & \ldots, & \alpha_{\ell(\alpha)},&1,&1 \\
           3,&\beta_3,  &\ldots, &\beta_{\ell(\beta)}
         \end{array}\right)_3.
\end{equation}
Clearly, $d'=3$, $\delta_1=3=d'$, $\gamma_1=\alpha_4\le 2<d'$ and $\ell(\gamma)-\ell(\delta)=(\ell(\alpha)-1)-(\ell(\beta)-1)=0$. Moreover,
\begin{align*}
  |\kappa_2(\lambda)| &=2+\left(\sum_{i=4}^{\ell(\alpha)}\alpha_i\right)+3
  +\left(\sum_{i=3}^{\ell(\beta)}\beta_i\right)+9 \\
&=|\alpha|-6+2+|\beta|-4+3+9 \\
   &=|\alpha|+|\beta|+4=|\lambda|.
\end{align*}
Hence $\kappa_2(\lambda)\in U_5^2(0,n)$.

We proceed to show $\kappa_2$ is an injection. Let $K_5^2(0,n)$ denote the image set of $\kappa_2$, which has been   shown to be a subset of $U_5^2(0,n)$. Given $\mu=(\gamma,\delta)_{d'}\in K_5^2(0,n)$, from the construction of $\kappa_2$, we see that $d'=3$, $\delta_1=3>\delta_2$, $\gamma_1\le 2$ and $\gamma_{\ell(\gamma)}=\gamma_{\ell(\gamma)-1}=1$. Define the map $\zeta_5^2$ as follows:
\begin{align*}
 \zeta_5^2(\mu)  =\left(\begin{array}{c}
                      \alpha \\
                      \beta
                    \end{array}\right)_d
    =
\left(\begin{array}{cccccccccccc}
        2, & 2, & 2, & \gamma_1,&\ldots,&\gamma_{\ell(\gamma)-2} \\
        2,&2, &\delta_2&\ldots, &\delta_{\ell(\delta)} &
      \end{array}\right)_{2}.
\end{align*}
It is trivial to check that $\zeta_5^2(\mu)\in P_6(0,n)$ and for any $\lambda\in P_6(0,n)$, $\zeta_5^2(\kappa_2(\lambda))=\lambda$, which implies $\kappa_2$ is an injection.
\qed

For example, let
$$\lambda=\left(\begin{array}{ccccccccccc}
                  6, & 6, &6,& 6, & 5,&4,&3,&3,&1 \\
                  6,  &6,& 5, & 5, & 4,&4,&2,&1,&1
                \end{array}\right)_{6}\in P_5(0,110).$$
Applying $\kappa_1$ on $\lambda$, we see that $d\ge 3$, $t=9$ and $k=6$, thus
$$\kappa_1(\lambda)=
\left(\begin{array}{ccccccccccc}
                  6,& 5, & 5, & 4,&4,&2,&1,&1 \\
                  7,& 6, & 5,&4,&4,&3,&3,&1
                \end{array}\right)_{7}$$
which is in $U_5^1(0,110)$. Applying $\zeta_5^1$ on $\kappa_1(\lambda)$, we recover $\lambda.$

For another example, let
\[\lambda=\left(\begin{array}{ccccccc}
                  2, & 2, & 2, & 2, & 2, & 1, & 1\\
                  2, & 2, & 2, & 2, & 1, & 1, & 1
                \end{array}\right)_2,\]
which is in $P_6(0,27)$. We see that
\[\kappa_2(\lambda)=\mu=\left(\begin{array}{ccccccc}
                      2, & 2, & 1, & 1,&1,&1 \\
                  3, &  2, & 2, & 1, & 1, & 1
                \end{array}\right)_3\]
which is in $U_5^2(0,27)$. Applying $\zeta_5^2$ on $\mu$, we recover $\lambda.$

We next present the injection $\chi_6$.

\begin{lem}\label{lem-4-8}
  For $n\ge 5$, there exists an injection $\chi_6$ from $P_8(0,n)$ into $U_4(0,n)$.
\end{lem}

\pf Given $\lambda=(\alpha,\beta)_d\in P_8(0,n)$. By definition, $\ell(\alpha)=\ell(\beta)$, $\beta_1<d$ and $\alpha_1\le d-2$. We claim that in this case $d\ge 3$. In fact if $d=2$, then $\alpha=\emptyset$. From $\ell(\alpha)=\ell(\beta)$ we see that $\beta=\emptyset$. Thus $|\lambda|=4$, contradicts to $n\ge 5$. Hence our claim has been verified.

Let $k=\max\{p\colon 1\le p\le \ell(\beta),\beta_p\ge 2\}$. If such $k$ not exists, then set $k=0$. Define
\begin{equation}\label{chi6}
  \chi_6(\lambda)=\left(\begin{array}{c}
                          \gamma \\
                          \delta
                        \end{array}\right)_{d'}
  =\left(\begin{array}{cccccccccccc}
           d-2,&\alpha_1, & \ldots, & \alpha_{\ell(\alpha)}&&& \\
           d-1,&\beta_1,  &\ldots, &\beta_{k},&2,&\beta_{k+1},&\ldots,&\beta_{\ell(\beta)}.
         \end{array}\right)_{d-1}.
\end{equation}
 It is routine to check that $\gamma$ and $\delta$ are partitions when one notice that $d-1\ge 2$. Moreover, we see that $\gamma_1=d-2=d'-1<d'$, $\delta_1=d-1=d'$ and $\ell(\gamma)-\ell(\delta)=(\ell(\alpha)+1)-(\ell(\beta)+2)=-1$. Furthermore,
 \begin{align*}
 |\chi_6(\lambda)| &=(d-2)+|\alpha|+(d-1)+2+|\beta|+(d-1)^2 \\
    &=|\alpha|+|\beta|+d^2=|\lambda|.
 \end{align*}
Hence $\chi_6(\lambda)\in U_4(0,n)$.

To show $\chi_6$ is an injection, let $K_6(0,n)$ denote the image set of $\chi_6$, which has been   shown to be a subset of $U_4(0,n)$. Given $\mu=(\gamma,\delta)_{d'}\in K_6(0,n)$, from the construction of $\chi_6$, we see that $d'=d-1\ge 2$,  $\gamma_1=d-2= d'-1$ and there exists $r\ge 2$ such that $\delta_r=2$. Define the map $\zeta_6$ as follows:
\begin{align*}
 \zeta_6(\mu)  =\left(\begin{array}{c}
                      \alpha \\
                      \beta
                    \end{array}\right)_d  =
\left(\begin{array}{cccccccccccc}
        \gamma_2, & \ldots, & \gamma_{\ell(\gamma)} &&& \\
        \delta_2,& \ldots, &\delta_{r-1},&\delta_{r+1}, &\ldots, &\delta_{\ell(\delta)}      \end{array}\right)_{d'+1}.
\end{align*}
It is trivial to check that $\zeta_6(\mu)\in P_8(0,n)$ and for any $\lambda\in P_8(0,n)$, we have $\zeta_6(\chi_6(\lambda))=\lambda$, which implies $\chi_6$ is an injection. \qed

For instance, let
$$\lambda=\left(\begin{array}{ccccccccccc}
                  4,&3,&3,&2,&2,&1,&1 \\
                   5, & 4, & 4,&3,&3,&2,&1
                \end{array}\right)_{6}\in P_8(0,74).$$
Applying $\chi_6$ on $\lambda$, we see that   $k=6$, thus
$$\chi_6(\lambda)=\mu=
\left(\begin{array}{ccccccccccc}
                 4,& 4,&3,&3,&2,&2,&1,&1 \\
                  5,&5, & 4, & 4,&3,&3,&2,&2,&1
                \end{array}\right)_{5}$$
which is in $U_4(0,74)$. Applying $\zeta_6$ on $\mu$, we recover $\lambda.$

The following lemma determines the injection $\chi_7$.

\begin{lem}\label{lem-chi7}
  For $n\ge 7$, there is an injection $\chi_7$ from $P_9(0,n)$ into $U_4(0,n)$.
\end{lem}

\pf For any $\lambda=(\alpha,\beta)_d\in P_9(0,n)$, by definition $\ell(\alpha)=\ell(\beta)$, $\beta_1<d$ and $\alpha_1=d-1>\alpha_2$. We first show that $d\ge 3$. If $d=1$, then  we see that $\alpha=\beta=\emptyset.$ Thus $|\lambda|=1$, contradicts to $n\ge 7$. Similarly, when $d=2$, we have $\alpha=\beta=(1)$, which implies $|\lambda|=6$, also contradicts to $n\ge 7$.  Hence $d\ge 3$. Define
\begin{equation}\label{chi7}
  \chi_7(\lambda)=\left(\begin{array}{c}
                          \gamma \\
                          \delta
                        \end{array}\right)_{d'}
  =\left(\begin{array}{cccccccccccc}
           d-2,&\alpha_1-1, & \alpha_2,&\ldots, & \alpha_{\ell(\alpha)},&1 \\
           d-1,&\beta_1,  &\ldots, &\beta_{\ell(\beta)},&1,&1
         \end{array}\right)_{d-1}.
\end{equation}
It is easy to see that $\gamma$ and $\delta$ are partitions. Moreover, $\delta_1=d-1=d'$, $\gamma_1=d-2=d'-1$ and $\ell(\gamma)-\ell(\delta)=\ell(\alpha)+2-(\ell(\beta)+3)=-1$.
Furthermore,
\begin{align*}
  |\chi_7(\lambda)| &= (d-2)+|\alpha|-1+1+(d-1)+|\beta|+2+(d-1)^2 \\
  &=|\alpha|+|\beta|+d^2=|\lambda|.
\end{align*}
Hence $\chi_7(\lambda)\in U_4(0,n)$. 

To show $\chi_7$ is an injection, let $K_7(0,n)$ denote the image set of $\chi_7$, which has been   shown to be a subset of $U_4(0,n)$. Given $\mu=(\gamma,\delta)_{d'}\in K_7(0,n)$, from the construction of $\chi_7$, we see that $d'=d-1\ge 2$, $\delta_1=d'$, $\gamma_1=\gamma_2=d'-1$. Moreover $s(\gamma)=1$ and $\delta_{\ell(\delta)}=\delta_{\ell(\delta)-1}=1$. Define the map $\zeta_7$ as follows:
\begin{align*}
 \zeta_7(\mu) =\left(\begin{array}{c}
                      \alpha \\
                      \beta
                    \end{array}\right)_d =
\left(\begin{array}{cccccccccccc}
        \gamma_2+1, &\gamma_3,& \ldots, & \gamma_{\ell(\gamma)-1} \\
        \delta_2,& \ldots, &\delta_{\ell(\delta)-2}  &   \end{array}\right)_{d'+1}.
\end{align*}
It is trivial to check that $\zeta_7(\mu)\in P_9(0,n)$ and for any $\lambda\in P_9(0,n)$, $\zeta_7(\chi_7(\lambda))=\lambda$, which implies $\chi_7$ is an injection. \qed

For instance, let
$$\lambda=\left(\begin{array}{ccccccccccc}
                  6,&5,&5,&4,&4,&2,&1 \\
                   4, & 4, & 4,&3,&3,&2,&1
                \end{array}\right)_{7}\in P_9(0,97).$$
Applying $\chi_7$ on $\lambda$,
$$\chi_7(\lambda)=\mu=
\left(\begin{array}{ccccccccccc}
                 5,& 5,&5,&5,&4,&4,&2,&1,&1 \\
                  6,&4, & 4, & 4,&3,&3,&2,&1,&1,&1
                \end{array}\right)_{6}$$
which is in $U_4(0,97)$. Applying $\zeta_7$ on $\mu$, we recover $\lambda.$

Finally, we   present the injection $\chi_8$ as follows.

\begin{lem}\label{lem-chi8}
For $n\ge 10$, there is an injection $\chi_8$ from $P_{10}(0,n)$ into $U_5(0,n)$.
\end{lem}

\pf Given $\lambda=(\alpha,\beta)_d\in P_{10}(0,n)$, by definition we have $\alpha_1=\alpha_2=d-1$, $\beta_1<d$ and $\ell(\alpha)=\ell(\beta)$. If $d=1$, we have $\alpha=\beta=\emptyset$. Thus $|\lambda|=1$,  contradicts to $n\ge 10$. Hence $d\ge 2$.  There are two cases.

Case 1: $d\ge 3$. In this case, set $\ell(\alpha)=\ell(\beta)=t$ and let $k=\max\{p\colon 1\le p\le t, \alpha_p\ge d-2\}$. From $\alpha_2=d-1$ we see that $k\ge 2$. Define
\begin{equation}\label{chi8-1}
  \chi_8(\lambda)=\left(\begin{array}{c}
                          \gamma \\
                          \delta
                        \end{array}\right)_{d'}
  =\left(\begin{array}{cccccccccccc}
            \beta_1,  &\ldots, &\beta_{t}&&&&& \\
          d,&\alpha_3,&\ldots, & \alpha_k,&d-2,&\alpha_{k+1},&\ldots, & \alpha_{t}.
         \end{array}\right)_{d}.
\end{equation}
From the choice of $k$, it is clear that $\delta$ is a partition. Moreover, $\gamma_1=\beta_1<d=d'$, $\delta_1=d=d'$ and $\ell(\gamma)-\ell(\delta)=t-t=0$. Furthermore, it is clear that $|\gamma|=|\beta|$ and $|\delta|=|\alpha|$. Hence $|\lambda|=|\chi_8(\lambda)|$. From the above analysis, we deduce that $\chi_8(\lambda)$ is in $U_5(0,n)$ with $d'\ge 3$.

Case 2: $d=2$. In this case, we see that $\beta_1=1$ and $\alpha_1=1$. Let $\ell(\alpha)=\ell(\beta)=t$, we have $\beta=\alpha=(1^t)$, which implies $n=2t+4$.
By $n\ge 10$ we have $t\ge 3$. Moreover, there is only one partition $\lambda=(1^t,1^t)_2$ lies in $P_{10}(0,n)$. Define
\begin{equation}\label{chi8-2}
  \chi_8(\lambda)=\left(\begin{array}{ccc}
                    1^{t-1} &  & \\
                    2, &2 ,& 1^{t-3}
                  \end{array}\right)_2.
\end{equation}
It is trivial to check that $\chi_8(\lambda)\in U_5(0,n)$ with $d'=2$.

We conclude this proof by show that $\chi_8$ is an injection. Since in Case 1 the image set satisfies $d'\ge 3$ and in Case 2 the image set satisfies $d'=2$, we see that the image set of the above two cases are not intersect. Note that Case 2 has only one partition,  it sufficient to show that $\chi_8$ in Case 1 is an injection. Let $K_8(0,n)$ denote the image set of $\chi_8$ in Case 1, which has been   shown to be a subset of $U_5(0,n)$. Given $\mu=(\gamma,\delta)_{d'}\in K_8(0,n)$, from the construction of $\chi_8$, we see that $d'\ge 3$, $\delta_1=d'>\delta_2$, and there exists $k$ such that $\delta_k=d'-2$. Define the map $\zeta_8$ as follows:
\begin{align*}
 \zeta_8(\mu)  =\left(\begin{array}{c}
                      \alpha \\
                      \beta
                    \end{array}\right)_d  =
\left(\begin{array}{cccccccccccc}
        d'-1, &d'-1,&\delta_2,& \ldots, & \delta_{k-1},&\delta_{k+1},&\ldots,&\delta_{\ell(\delta)} \\
        \gamma_1,& \ldots, &\gamma_{\ell(\gamma)}  & &&&&  \end{array}\right)_{d'}.
\end{align*}
It is trivial to check that $\zeta_8(\mu)\in P_{10}(0,n)$ and for any $\lambda\in P_{10}(0,n)$, $\zeta_8(\chi_8(\lambda))=\lambda$, which implies $\chi_8$ is an injection. \qed

For example, let
$$\lambda=\left(\begin{array}{ccccccccccc}
                  4,&4,&4,&4,&3,&2,&2 \\
                   4, & 4, & 3,&2,&2,&2,&1
                \end{array}\right)_{5}\in P_{10}(0,66).$$
Applying $\chi_8$ on $\lambda$, we see that $d\ge 3$ and $k=5$. Thus
$$\chi_8(\lambda)=\mu=
\left(\begin{array}{ccccccccccc}
                 4, & 4, & 3,&2,&2,&2,&1 \\
                  5,&4,&4,&3,&3,&2,&2
                \end{array}\right)_{5}$$
which is in $U_5(0,66)$. Applying $\zeta_8$ on $\mu$, we recover $\lambda.$

Now we are in a position to prove Theorem \ref{thm-est-u0n-u1n} with the aid of Lemma \ref{lem-chi2}--Lemma \ref{lem-chi8}.

{\noindent \it Proof of Theorem \ref{thm-est-u0n-u1n}.} It is trivial to check that for $6\le n\le 14$, Theorem \ref{thm-est-u0n-u1n} holds. Recall that in the proof of Lemma \ref{lem-4-8}, we build an injection $\chi_6$ from $P_8(0,n)$ into $U_4(0,n)$. The image set of $\chi_6$ is denoted by $K_6(0,n)$. We next show that for $n\ge 15$,   there exists $\mu\in U_4(0,n)$ and not in $K_6(0,n)$, which implies for $n\ge 15$,
\begin{equation}\label{ine-str-ospt}
  \# U_4(0,n)>\# P_8(0,n).
\end{equation}
In fact, if $n=2t$ for some $t\ge 8$, define
\begin{equation*}
  \mu=\left(\begin{array}{c}
              \gamma \\
              \delta
            \end{array}\right)_{d'}=\left(\begin{array}{ccccc}
              1^{t-6}  \\
              3, & 1^{t-6}
            \end{array}\right)_3.
\end{equation*}
Clearly, $\mu\in U_4(0,n)$. Moreover, from $\gamma_1=d'-2$ we see that $\mu\not\in K_6(0,n)$.
If  $n=2t+1$ for some $t\ge 7$, then define
\begin{equation*}
  \mu=\left(\begin{array}{c}
              \gamma \\
              \delta
            \end{array}\right)_{d'}=\left(\begin{array}{ccccc}
              1^{t-6}  \\
              3, & 2,&1^{t-7}
            \end{array}\right)_3.
\end{equation*}
Similarly, $\mu\in U_4(0,n)$ but by $\gamma_1=d'-2$, we have $\mu\not\in K_6(0,n)$. Thus \eqref{ine-str-ospt} holds.

From Lemma \ref{lem-chi2}--Lemma \ref{lem-chi8} and \eqref{ine-str-ospt}, we see that
\begin{align*}
  \#P_1(0,n) &= \# U_5(0,n) \\
 \#P_2(0,n) &= \# U_5(0,n)\\
   \#P_3(0,n) &\le \# U_4(0,n) \\
  \#P_4(0,n)+\#P_7(0,n) &\le \# U_4(0,n)\\
 \#P_5(0,n)+\#P_6(0,n) &\le \# U_5(0,n)\\
  \#P_8(0,n) &< \# U_4(0,n) \\
 \#P_9(0,n) &\le \# U_4(0,n)\\
 \#P_{10}(0,n) &\le \# U_5(0,n)
\end{align*}
Thus by the above equations and inequalities, we deduce \eqref{ine-u4n-u5n-n0n}. This completes the proof. \qed

\section{Concluding Remarks}\label{6}

We end this paper by rasing the following conjecture.

\begin{conj}\label{conj-ospt}
Can one prove that $u(0,n)-u(1,n)\ge \frac{N(0,n)}{2}$ for $n\ge 8$? This inequality has been verified for $8\le n\le 10000$.
\end{conj}

It should be remarked that by Corollary \ref{cor-u0n-ospt}, Conjecture \ref{conj-ospt}  implies
      \[\ospt(n)\ge \frac{p(n)}{4}+\frac{N(0,n)}{2},\]
which improves both \eqref{ospt-lower} and \eqref{chan-mao-lowerbound}. Moreover, by \eqref{ospt-ays}, the coefficient $1/2$ of $N(0,n)$ in Conjecture \ref{conj-ospt} is sharp.

\vskip 0.5cm

\noindent{\bf Acknowledgments.} The  authors thank Nianhong Zhou for his enlightening discussions, from which this work has benefited a lot. The  author was supported by the National Science Foundation of China grants 12171358 and 12371336.

\end{document}